\documentclass{amsart}\usepackage{amsfonts,amssymb,bbm}

\def\k{\mathbbm{k}}
\def\N{\Bbb{N}}
\def\Q{\Bbb{Q}}\def\R{\Bbb{R}}\def\Z{\Bbb{Z}}

\def\li{\ \\ $\bullet$ }\newcommand{\ls}{~\\ $\star$ }\def\di{\partial}

\def\suml{\sum\limits}
\def\capl{\mathop\cap\limits}\def\cupl{\mathop\cup\limits}

\newcommand{\quotients}[2]{{\footnotesize\left.\raisebox{0.4ex}{$#1$}\! / \!\raisebox{-0.4ex}{$#2$}\right.}}

\def\hF{\hat{F}}\def\hG{{\hat{G}}}

\def\hW{{\widehat{W}}}\def\hy{\hat{y}}

\def\de{\delta}\def\De{\Delta}

\def\la{\lambda}

\def\cA{\mathcal A}\def\ca{\frak a}
\def\cD{\mathcal D}
\def\cG{\mathcal G}\def\cK{{\mathcal K}}\def\cR{\mathcal R}
\def\cm{{\frak m\hspace{0.05cm}}}

\def\uu{{u}}\def\uv{{v}}\def\uy{{y}}

\def\one{{1\hspace{-0.1cm}\rm I}}

\def\bx{{x}}\def\by{{y}}\def\bz{{z}}

\newcommand{\ber}{\begin{array}{l}}\newcommand{\eer}{\end{array}}
\newcommand{\bpm}{\begin{pmatrix}}\newcommand{\epm}{\end{pmatrix}}
\newcommand{\bM}{\begin{matrix}}\newcommand{\eM}{\end{matrix}}
\newcommand{\bee}{\begin{enumerate}}\newcommand{\eee}{\end{enumerate}}

\def\wrt{with respect to }\def\sset{\subset}\def\sseteq{\subseteq}\def\ssetneq{\subsetneq}

\def\Mat{Mat(m,n;R)}

\def\obo{order-by-order}\def\wt{good}

\newcommand{\beq}{\begin{equation}}\newcommand{\eeq}{\end{equation}}
\newtheorem{Lemma}{Lemma}[section]\newcommand{\bel}{\begin{Lemma}}\newcommand{\eel}{\end{Lemma}}
\newtheorem{Example}[Lemma]{Example}\newcommand{\bex}{\begin{Example}\rm}\newcommand{\eex}{\end{Example}}
\newtheorem{Proposition}[Lemma]{Proposition}\newcommand{\bprop}{\begin{Proposition}}\newcommand{\eprop}{\end{Proposition}}
\newtheorem{Definition-Proposition}[Lemma]{Definition-Proposition}

\def\bpr{~\\{\em Proof.\ }}
\newcommand{\epr}{$\blacksquare$\\}
\newtheorem{Theorem}[Lemma]{Theorem}\newcommand{\bthe}{\begin{Theorem}}\newcommand{\ethe}{\end{Theorem}}
\newtheorem{Definition}[Lemma]{Definition}\newcommand{\bed}{\begin{Definition}}\newcommand{\eed}{\end{Definition}}
\newtheorem{Remark}[Lemma]{Remark}\newcommand{\beR}{\begin{Remark}\rm}\newcommand{\eeR}{\end{Remark}}
\newtheorem{Corollary}[Lemma]{Corollary}\newcommand{\bcor}{\begin{Corollary}}\newcommand{\ecor}{\end{Corollary}}
\newcommand{\bet}{\begin{tabular}{cccccccc}}\newcommand{\eet}{\end{tabular}}

\newcommand{\bin}[2]{\binom{#1}{#2}}
\def\liml{\lim\limits}

\title[]{A\MakeLowercase{ strong version of implicit function theorem}}
\author[]{G\MakeLowercase{enrich} B\MakeLowercase{elitskii and} D\MakeLowercase{mitry} K\MakeLowercase{erner}}
\address{Department of Mathematics, Ben Gurion University of the Negev, P.O.B. 653, Be'er Sheva 84105, Israel.}
\email{genrich@math.bgu.ac.il, dmitry.kerner@gmail.com}
\date{\today}
\thanks{D.K. was partially supported by the grant FP7-People-MCA-CIG, 334347 and the grant ISF 844/14.}

\begin{document}
\begin{abstract}
We suggest the {\em necessary/sufficient} criteria for the
existence of a (\obo) solution $\by(\bx)$ of a functional equation $F(\bx,\by)=0$ over a  ring. In full generality, the criteria
hold in the category of filtered groups, this includes the wide class of modules over (commutative, associative) rings.
The classical implicit function theorem and its strengthening obtained by Tougeron and Fisher
appear to be (weaker) particular forms of the general criterion.

We obtain a special criterion for solvability of the equations arising from group actions, $g(w)=w+u$, here $u$ is ``small". As an immediate application we re-derive the classical criteria of determinacy, in terms of the tangent space to the orbit.

Finally, we prove the Artin-Tougeron-type approximation theorem: if a system of $C^\infty$-equations has a formal solution and the derivative satisfies a Lojasiewicz-type condition then the system has a $C^\infty$-solution.
\end{abstract}
\maketitle
\setcounter{secnumdepth}{6} \setcounter{tocdepth}{1}\vspace{-0.25cm}
\tableofcontents

All the rings in this paper are commutative, associative, with unit element, of zero characteristic.

\section{Introduction}
We use the multivariable notation, $\bx=(x_1,\dots,x_m)$, $\by=(y_1,\dots,y_n)$.
\subsection{The general setting and known results}

Consider a system of (analytic/formal/$C^\infty$/$C^k$) equations $F(\bx,\by)=0$.
The classical Implicit Function Theorem reads:
\begin{equation*}\ber
\text{if the matrix of derivatives, $\di_\by F(0,0)=F'_\by(0,0)$, is right
  invertible (i.e. is of the full rank)}
\\\text{then $F(\bx,\by)=0$, has a (analytic/formal/etc.) solution.}
\eer\end{equation*}
The condition ``$F'_\by(0,0)$ is right invertible" is quite restrictive. For example, the
 theorem does not ensure a solution of the one-variable equation $xy=0$ (in the vicinity of $(0,0)$) or of $y^2=0$ (at any point).

Various strengthenings/generalizations of this theorem are known (including Hensel lemma). For example,
Tougeron's implicit function theorem ensures solvability when the matrix $F'_{\by}(\bx,0)$ is not too degenerate.
Denote by $I_{max}F'_{\by}(\bx,0)$ the ideal of the maximal minors of this matrix.
\bthe\cite[page 56]{Tougeron-book},\cite{Tougeron1968}
Let $R=\k[[\bx,\by]]$ or $\k\{\bx,\by\}$ (for $\k$ a normed field) or $C^\infty(\R^m\times \R^n,0)$.
Let $F(\bx,\by)\in R^{\oplus p}$, with $p\le n$, and let $I\sset R$ be a proper ideal.
 If $F(\bx,0)\in I\cdot\Big(I_{max}F'_{\by}(\bx,0)\Big)^2R^{\oplus p}$
 then there exists a solution, $F(\bx,\by(\bx))\equiv0$,
 satisfying the condition: $\by(\bx)\in I R^{\oplus n}$.
\ethe

While this theorem ensures the solution of $yx=0$ and $y^2=0$, it fails to ensure the solution of the
 system
\begin{equation*}
 \Big\{\ber y^2_1+y_1x=x^3\\y^2_2+y_2x=x^3\eer.
\end{equation*}
Here $F(x,0)=x^3\bin{1}{1}$, $I_{max}F'_\by(x,0)=(x^2)$, thus $F(x,0)\not\in \Big(I_{max}F'_\by(x,0)\Big)^2$.

Tougeron himself realized in \cite{Tougeron1966} that one can replace
 in the condition $F(\bx,0)\in I\cdot\Big(I_{max}F'_{\by}(\bx,0)\Big)^2R^{\oplus p}$ the ideal $I_{max}F'_{\by}(\bx,0)$ by
  the larger ideal,
\begin{equation*}
\ca_{F'_{\by}(\bx,0)}:=ann\Big(coker(F'_{\by}(\bx,0))\Big),
\end{equation*}
 the {\em annihilator of
  the cokernel of the morphism} $R^{\oplus n}\stackrel{F'_{\by}(\bx,0)}{\longrightarrow}R^{\oplus p}$.
  Some properties of this ideal are given in \S\ref{Sec.Background.Ann.Coker.and.J}. By now we just mention
  that for $p=1$, i.e. the case of one equation, the two ideals coincide: $\ca_{F'_{\by}(\bx,0)}=I_{max}F'_{\by}(\bx,0)$.

The statement was further strengthened by B.Fisher, he replaced one of the factors in
 $\Big(\ca_{F'_{\by}(\bx,0)}\Big)^2$ by the image $Im(F'_\by(\bx,0))\sseteq R^{\oplus p}$.
(The initial version was for p-adic rings, we give a more general version relevant to our context.)
\bthe\cite{Fisher1997}\label{Thm.Fisher}
Let $(R,\cm)$ be a local Henselian ring over a field of zero characteristic. Let $F_1,\dots,F_p\in R[[y_1,\dots,y_p]]$.
 Suppose
\beq\label{Eq.Fisher's.Condition}
F(\bx,0)\in \cm\cdot \ca_{F'_\by(\bx,0)}\cdot Im\Big(F'_\by(\bx,0)\Big).
\end{equation}
  Then there exists a solution, $F(\bx,\by(\bx))\equiv0$, satisfying $\by(\bx)\in\cm\cdot \ca_{F'_{\by}(\bx,0))}\cdot R^{\oplus p}$.
\ethe

In the case of one equation,  $p=1$, this coincides with Tougeron's result. For $p>1$ Fisher's result is stronger.
(Note that $\ca_{F'_{\by}(\bx,0)}Im\Big(F'_\by(\bx,0)\Big)\supseteq \Big(\ca_{F'_\by(\bx,0)}\Big)^2 R^{\oplus p}$,  and for $p>1$ the inclusion  is in general proper.)

Though Fisher's version solves the examples mentioned above, it cannot cope with a slightly more complicated example:
\beq\label{Eq.hard.equation}
 y^2_1-y^2_2+y_1x^k_1+y_2x^k_2+g(x_1,x_2)=0,
\eeq
where  $\cm=(x_1,x_2)\sset R=\k[[x_1,x_2]]$ and $g(x_1,x_2)\in \cm^{2k+1}$ for $k>2$.
Here
\begin{equation*}
\cm\cdot\ca_{F'_\by(\bx,0)}\cdot Im\Big(F'_\by(\bx,0)\Big)=\cm\cdot(x^k_1,x^k_2)^2\subsetneq\cm^{2k+1},
\end{equation*}
 thus in general $F(\bx,0)=g(x_1,x_2)\not\in\cm\cdot\ca_{F'_\by(\bx,0)}\cdot Im\Big(F'_\by(\bx,0)\Big)$, i.e. the condition \eqref{Eq.Fisher's.Condition} is not satisfied.

\subsection{Overview of the results}\label{Sec.Intro.Sketch.of.Results}
Our work has began from the observation that Fisher's condition can be further weakened: instead of
  $F(\bx,0)\in \cm\cdot \ca_{F'_\by(\bx,0)}\cdot Im\Big(F'_\by(\bx,0)\Big)$ it is enough to ask for
 $F(\bx,0)\in \cm\cdot J\cdot Im\Big(F'_\by(\bx,0)\Big)$, where $\ca_{F'_\by(\bx,0)}\sseteq J\sset R$ is the biggest possible ideal satisfying
$J^2=J\cdot\ca_{F'_\by(\bx,0)}$. (See corollary \ref{Thm.IFT.J^2=JA} and \S\ref{Sec.Comparison.to.Fisher.Tougeron}  for more detail.)
This gives the further strengthening of Fisher's and Tougeron's
 statements. Still, this strengthening does not help to address the very simple system (cf. \S\ref{Sec.Comparison.to.Fisher.Tougeron})
\begin{equation*}
\Big\{\ber y^2_1+y_1x_1=x^{100}_1\\y^2_2+y_2x_2=x^{100}_2\eer.
\end{equation*}
Here $\ca_{F'_\by}=(x_1x_2)$ thus $J=(x_1x_2)$, but $F(\bx,0)\not\in J\cdot Im\Big(F'_\by(\bx,0)\Big)$.

\

In this note we prove much stronger solvability criteria. Here we sketch just the main features of the method.
The detailed formulation can be found in \S\ref{Sec.Criteria.General} (theorem \ref{Thm.IFT.Nec.Suff.Criterion.non.Abelian.Groups}) and \S\ref{Sec.Approximation.Theorem} (theorem \ref{Thm.Approximation.Theorem.for.smooth.funcs}), the applications are in
 \S\ref{Sec.Remarks} and \S\ref{Sec.Group.Action.Equations.applications}.

\vspace{-0.25cm}
\subsubsection{} We weaken the condition on $F'_\by(\bx,0)$ further, to the ``weakest possible" condition of ``iff" type, so that we get a Strong Implicit Function Theorem.

 Our results hold in broader category.
  It is natural to extend from the classical case of $\k[[\bx,\by]]$, $\k\{\bx,\by\}$, $C^p(\R^m_\bx\times\R^n_\by,0)$ to the local Henselian rings
 (not necessarily regular or Noetherian) over a field. In fact even the ring structure is not necessary, our main result, theorem \ref{Thm.IFT.Nec.Suff.Criterion.non.Abelian.Groups}, is for the filtered (not necessarily abelian) groups.\vspace{-0.25cm}
\subsubsection{} A particular class of equations comes from the group actions, $G\circlearrowright W$. Assume $W$ is a filtered abelian group (e.g. a module over a local ring). To understand how large is the orbit one studies the equation  $g(w)=w+u$. Here $g\in G$ is an unknown, while $u\in W$ is ``small". (More precisely, one studies whether the orbit $Gw$ is open in the topology defined by the filtration.)
Theorem  \ref{Thm.IFT.Nec.Suff.Criterion.non.Abelian.Groups}, being very general, is of little use here. Rather, we obtain a special version of strong IFT, \S\ref{Sec.Group.Action.IFT}.
\vspace{-0.25cm}
\subsubsection{}  Usually the main problem is to establish  the \obo\ solution procedure. Thus many of our results are of the form
 ``If (\dots) then there exists a Cauchy sequence $\{\by^{(n)}(\bx)\}_n$ such that $F(\bx,\by^{(n)}(\bx))\to0$".
(The topology here is induced by a filtration, e.g. by the powers of maximal ideal.)

Once such a result is established, one has a solution in the completion of $R^{\oplus p}$ by the filtration. Then (if $R$ is non-complete) one uses the
 Artin-type approximation theorems, \cite{K.P.P.R.M.}, to establish a solution over $R$, or at least over the henselization of $R$.

For the ring $C^\infty(\R^p,0)$ and many other important rings the Artin approximation does not hold (in the naive way).
  Over some rings we can directly ensure a solution, see \S\ref{Sec.Exact.Solutions.Over.Rings.with.IFT}.
 For $C^\infty(\R^p,0)$ we use theorem \ref{Thm.Approximation.Theorem.for.smooth.funcs}.

\

\subsection{Comments and motivation}
Several remarks/explanations are necessary at this point.

 \subsubsection{}
  Recall the simple geometric interpretation. Consider the (germ of) subscheme/subspace\\ $\{F(\bx,\by)=0\}\sset(X,0)\times(Y,0)$. The classical IFT,  in the case when $(X,0)$, $(Y,0)$ are smooth, gives a sufficient condition that the germ $\{F(\bx,\by)=0\}$ is smooth and its projection onto $(X,0)$ is an isomorphism. Our version of IFT,
 for the arbitrary henselian germs $(X,0)$, $(Y,0)$, gives a {\em necessary and sufficient} condition that
the germ $\{F(\bx,\by)=0\}$ has an irreducible component whose projection onto $(X,0)$ is an isomorphism.

 This can be restated as follows. Consider the natural projection $\{F(\bx,\by)=0\}\stackrel{\pi}{\to}(X,0)$. Usually  this  projection is not an isomorphism. The solvability of the equation means the weaker property: the existence of the section of $\pi$, $(X,0)\stackrel{s}{\to}\{F(\bx,\by)=0\}$.


To emphasize: as the germ  $\{F(\bx,\by)=0\}$ is in general non-smooth (possibly reducible, non-reduced), the question cannot be simply ``linearized" by an automorphism of $(X,0)\times(Y,0)$, i.e. cannot reduced to the classical IFT by some appropriate change of variables.

\subsubsection{}
 A reformulation in terms of commutative algebra.
 Given a ring $R$ over some base ring $R_X$, e.g. $R=R_X[\by]$ or $R=R_X[[\by]]$, etc.
 Given an ideal $F=(F_1,\dots,F_p)\sset R$, a solution of $F(\bx,\by)=0$ is a projection
$R\stackrel{\by\to\by(\bx)}{\to}R_X$ whose kernel is precisely $F$.

\subsubsection{} The classical approach to construct a solution is the \obo\ approximation: first solve the part linear in $\by$ (modulo quadratic terms),
 then quadratic, cubic etc. Accordingly we always present the equation(s) $F(\bx,\by)=0$ in the form $\uu+L\by+H(\by)=0\in W$. Here  $\uu=F(\bx,0)\in W$,

 $V\stackrel{L=F'_\by(\bx,0)}{\longrightarrow}W$ is a homomorphism of $R$-modules
 (or just of abelian groups);

$H(\by)$ denotes the remaining ``higher order terms" (a contractive map in the sense of Krull topology), defined in \S\ref{Sec.Background.Iimpl.Func.Equation}.
\\Further, as we always start from a solution of the linear part, $\uu+L\by=0$, we assume $\uu\in L(V)$,
i.e. $\uu=-L\uv$, for some $\uv\in V$. Therefore the equation to solve is presented in the form
\begin{equation*}
L(\by-\uv)+H(\by)=0.
\end{equation*}
\subsubsection{} In practice one usually needs not just {\em a solution}. Thinking of $\uv$ as a parameter, one needs a statement of the type:
\begin{equation*}\ber
\text{There exists a subgroup/submodule $V_1\sseteq V$ such that for any $\uv\in V_1$
 the equation $L(\by-\uv)+H(\by)=0$ has}\\\text{a solution, $\by_\uv\in V_1$ which is ``close" to $\uv$
 and depends on $\uv$ ``differentiably".}\vspace{-0.2cm}
\eer\end{equation*}\vspace{-0.1cm}
 We call this a {\em \wt\ solution},  the precise formulation is in \S\ref{Sec.Background.Iimpl.Func.Equation}.
 Our criteria answer the question:
\begin{equation*}
\text{Given $V\stackrel{L,H}{\longrightarrow}W$, what is the biggest $V_1\sseteq V$ such that for any $v\in V_1$ there exists a good solution?}
\end{equation*}
Note that for some equations all the solutions are ``not good", cf. \S\ref{Sec.Equations.with.no.good.solutions}.

\subsubsection{} If the number of unknowns equals the number of equations and $F'_\by(0,0)$ is non-degenerate, then the classical IFT ensures the unique solution. When $F'_\by(0,0)$ is degenerate the solution (if it exists) can be non-unique, as the space $\{F(\bx,\by)=0\}$ can have several irreducible components. However, when $L$ is injective, the solution lying in $V_1$ is unique! The (non-)uniqueness issues are addressed in \S\ref{Sec.Uniqueness}.
\vspace{-0.25cm}
\subsubsection{} We expand $F(\bx,\by)=0$ in powers of $\by$  (i.e. at the point $\by=0$), hence the criteria are formulated in terms of $F(\bx,0)$,
 $F'_\by(\bx,0)$ etc. One can expand at some other point, $\by=\by^{(0)}(\bx)$, then the criteria are written in terms of
 $F(\bx,\by^{(0)}(\bx))$, $F'_{\by^{(0)}(\bx)}(\bx,\by^{(0)}(\bx))$,\dots  (For example, theorem \ref{Thm.Fisher} is stated in \cite{Fisher1997}
 in such a form.) Such an expansion at $\by^{(0)}(\bx)$ is helpful if one has a good initial approximation for the solution. The two approaches
 are obviously equivalent, e.g. by changing the variable $\by\to\by-\by^{(0)}(\bx)$.
  To avoid cumbersome formulas we always expand at $\by=0$.
\subsubsection{}
In view of our initial result, \S\ref{Sec.Intro.Sketch.of.Results}, one might try to weaken the condition on the ideal $\ca_{F'_\by(\bx,0)}\sseteq J\sset R$ further. It appears that $J^2=J\ca_{F'_\by(\bx,0)}$ is almost the ``weakest possible" among the conditions stated in terms of ideals only, it cannot be significantly weakened, cf. \S\ref{Sec.Comparison.of.ideal.conditions.to.filtration}. But this condition is still far from being necessary. The ``right" condition (necessary and sufficient) is obtained by replacing the ideals with filtered subgroups. As a bonus we do not need the rings structure anymore, e.g. theorem \ref{Thm.IFT.Nec.Suff.Criterion.non.Abelian.Groups} holds in the generality of (non necessarily abelian) filtered groups.
\subsubsection{} If the equations $F(\bx,\by)=0$ are linear in $\by$, i.e. $F(\bx,\by)=F(\bx,0)+F'_{\by}(\bx,0)\by$,
 then the obvious sufficient condition for solvability is: the entries of $F(\bx,0)$ lie
  in the ideal $\ca_{F'_{\by}(\bx,0)}$. While the (tautological) necessary and sufficient condition
  is: $F(\bx,0)\in Im(F'_{\by}(\bx,0))$. This condition is much weaker than those of Tougeron and Fisher and is far from being sufficient for non-linear equations.
  Therefore as landmarks for our criteria one should consider
  equations that are {\em non-linear in $\by$.}
\subsubsection{}
Implicit Function Theorem is a fundamental result. In \S\ref{Sec.Applications.Bifurcations} we obtain an immediate corollary to non-bifurcation of multiple polynomial roots  under deformations.
In \S\ref{Sec.Applications.Smooth.Germs} we indicate a potential application to the study of smooth curve-germs (lines/arcs) on singular spaces.
In \S\ref{Sec.Group.Action.Determinacy.Appplications} we apply a version of strong IFT to group-actions to re-derive the classical criteria of finite determinacy.

The further directions in Algebra and Geometry are:
 matrix equations, equations on (filtered) groups \cite{Bandman-Garion-Kunyavskii}, \ tactile maps  \ \cite{Bruschek-Hauser}, bounds on Artin-Greenberg functions \cite{Rond.2005}, \cite{Rond.2010}, etc.
We hope to report on these applications soon.

\subsection{Acknowledgements}
We thank J. Bochnak, H. Hauser, D. Popescu, J.M. Ruiz, E. Shustin, S. Yakovenko for the attention and valuable suggestions. We also thank the two referees, their numerous remarks have greatly improved the exposition.

\section{Definitions and Notations}

\subsection{Groups with descending filtration}\label{Sec.Background.Filtered.Groups}
We always assume that a (not necessarily abelian) group $V$ is filtered  by a sequence of normal subgroups $V\supset V_1\supset V_2\supset\cdots$, $V_j\lhd V$. Moreover, we assume that the filtration satisfies: $[V_1,V_i]\sseteq  V_{i+1}$, similarly to the lower central series of a group. This later condition is trivial when $V$ is an abelian group.
If $V$ is complete \wrt\ $\{V_i\}$ then
the filtration is faithful, i.e. $\capl_{i\ge1}V_i=\{{\bf 1}_V\}$.
The filtration induces the Krull topology, the fundamental system of neighborhoods of $v\in V$ is $\{v V_j\}_{j\ge1}$, or $\{V_jv\}_{j\ge1}$, by the normality.
\bex
1. The simplest case is when $V$ is a module over a ring, with filtration defined by the powers of an ideal, $V_j=I^j V$.
\\2. Let $(R,\cm)$ be a local ring with the filtration $R\supset\cm\supset\cm^2\supset\cdots$.
 Consider the group of invertible $R$-matrices, $V=GL(n,R)$. We get the filtration by the normal subgroups $V_j:=\{\one+A|\ A\in Mat(n,n;\cm^j)\}$.
\\3. Let $(X,0)$ be the germ of a space (algebraic/formal/analytic etc). Consider the group of its automorphisms, $V=Aut(X,0)$. The natural filtration is by the subgroups of automorphisms that are identity up to $j$'th order. More precisely, denote by $(R_{(X,0)},\cm)$ the local ring of (germs of) regular functions. Then
\begin{equation*}
Aut^{(j)}(X,0)=\{\phi\circlearrowright R_{(X,0)}|\ \phi(f)-f\in\cm^{j+1},\ \forall\ f\in R\}.
\end{equation*}
\eex

\subsection{Implicit function equation}\label{Sec.Background.Iimpl.Func.Equation}
Given two (not necessarily abelian) groups, $V$, $W$, a homomorphism $V\stackrel{L}{\to}W$, and a decreasing filtration $\{V_j\}$ by normal subgroups, we define the filtration $W_j:=L(V_j)$ on $W$.
\label{Sec.Background.Higher.Order.Terms} Consider the equation, $L(\uy)H(\uy)=L(\uv)$, where the ``higher order" map $V\stackrel{H}{\to} W$ , usually not a homomorphism, satisfies:
\li  $H({\bf 1}_V)={\bf 1}_W$ and
\li $\big(H(y)\big)^{-1}H(yV_j)\sseteq L(V_{j+1})$ for any $y\in V_1$ and any $j\in\N$.

Note that being of higher order depends essentially on $L$, in particular $H(V_j)\sseteq L(V_{j+1})$.

\subsubsection{}
If $V,W$ are abelian groups, then the implicit function equation is $L(\uy-\uv)+H(\uy)=0$, where $L\in Hom(V,W)$, while the higher order $H(\uy)$ satisfies: $H(0_V)=0_W$ and $H(y+V_j)-H(y)\sseteq L(V_{j+1})$ for any $y\in V_1$ and $j$.

Most common case is when $V,W$ are modules over a (commutative, associative) ring $R$.
Then usually $L\in Hom_R(V,W)$.
We say that the map $V\stackrel{H}{\to}W$ is of order $\ge k$ if for any ideal $J\sset R$ holds $H(JV)\sseteq J^k W$.
\bex
Suppose $R$ is graded, fix an ideal $J\sset R$, and consider the filtration $V_j=J^j V$. Suppose $H(y)$ can be written as a sum of homogeneous forms, $H(y)=\suml_{i\ge k} h_i(y)$, the degree of $h_i(y)$ being $i$. If $L(V)\supseteq J^{k-2}W$ then $H(y)$ is a ``higher order" term for $L$. Indeed, for any $i\ge k$ and $y\in J V$:
\begin{equation*}
h_i(y+V_j)-h_i(y)\in J^{i-1}\cdot J^j W\sseteq J^{k-1+j}W\sseteq J^{j+1}L(V).
\end{equation*}
\eex
\bex (Warning) Being of higher order terms can be  a restrictive condition. For example, in the equation $y^2-yx+x^a=0$ the monomial $y^2$ represents the higher order term for the filtration $V_j=(x^j)$ only if $a\ge3$. Otherwise the condition $H(V_1)\sseteq L(V_2)$ is not satisfied.
\eex

\subsubsection{}\label{Sec.Background.Order.By.Order.Solution}
An {\em \obo\ solution} of the equation $L(\uy)H(\uy)=L(\uv)$ is a Cauchy sequence, $\{\uy^{(n)}\}_{n\ge1}$ \wrt the filtration $V_\bullet$, i.e. $\uy^{(n)}(\uy^{(n+1)})^{-1}\in V_n$, such that  $L(\uy^{(n)})H(\uy^{(n)})L(\uv)^{-1}\in L(V_n)$.
By the normality, $V_n\lhd V$, we can also write the condition as $(\uy^{(n+1)})^{-1}\uy^{(n)}\in V_n$ or
$L(\uv)^{-1}L(\uy^{(n)})H(\uy^{(n)})\in L(V_n)$.

\subsubsection{}\label{Sec.Background.Good.Solution}  We say that the equation $L(\uy)H(\uy)=L(\uv)$ admits a {\em ``\wt\ solution on $V_1$"} if there exists a map $V_1\stackrel{\uy}{\to}V_1$ satisfying (we denote $\uy(v)$ by $\uy_v$, i.e. consider $v$ as a parameter):
\ls $L(\uy_\uv)H(\uy_\uv)=L(\uv)$ for any $\uv\in V_1$;
\ls  $\uy_{{\bf 1}_V}={\bf 1}_V$ and $y$ respects the filtration, $y(V_i)\sseteq V_i$, (this is a strengthening of continuity);
\ls the map $\uy$ is ``differentiable and close to identity", namely $y_v=g(v)\cdot v$ where the map $V_1\stackrel{g}{\to}V_1$ satisfies: $g(v\cdot V_j)g^{-1}(v)\in V_{j+1}$ for any $v\in V_1$ and $j\in\N$. Alternatively this condition can be stated as: $y_{v\De_j}\De^{-1}_jy_v^{-1}\in V_{j+1}$ for any $\De_j\in V_j$. By the normality this is equivalent to $\De^{-1}_jy_v^{-1}y_{v\De_j}\in V_{j+1}$.

\

We say that a solution $V_1\stackrel{y}{\to}V_1$ is {\em quasi-good} if $y_v=g(v)\cdot v$, where $g(V_j)\sseteq V_{j+1}$. (Good implies quasi-good.)

\subsubsection{}
Combining these notions we get the notion of a {\em\wt\ \obo}\ solution: a Cauchy sequence of maps, $\{V_1\stackrel{\uy^{(n)}}{\to}V_1\}_n$, with
 \li $\uy^{(n)}_{{\bf 1}_V}={\bf 1}_V$, $L(\uy^{(n)}_v)H(\uy^{(n)}_v)L(\uv)^{-1}\in L(V_n)$ and
\li $\forall n,j\ge1,\ \forall \De_j\in V_j$: \ $y^{(n)}_{v\De_j}\De^{-1}_j(y^{(n)}_v)^{-1}\in V_{j+1}$.

Similarly, a {\em quasi-good \obo\ solution} satisfies $y^{(n)}_v=g^{(n)}(v)\cdot v$, where $g^{(n)}(V_j)\sseteq V_{j+1}$.

\

If $V$, $W$ are abelian groups then all the notions simplify accordingly. A \wt\ \obo\  solution means a Cauchy sequence of maps $\{V_1\stackrel{\uy^{(n)}}{\to}V_1\}_n$ satisfying the conditions
\li
$\uy^{(n)}_{\uv}-\uy^{(n+1)}_{\uv}\in V_n$,\quad $L(\uy^{(n)}_{\uv}-\uv)+H(\uy^{(n)}_{\uv})\in L(V_n)$ and
\li $\forall n,j\ge1,\ \forall \De_j\in V_j:\ \uy^{(n)}_{\uv+\De_j}-\uy^{(n)}_{\uv}-\De_j\in V_{j+1}$.

\subsection{Annihilator of cokernel}\label{Sec.Background.Ann.Coker.and.J}
Consider a homomorphism of finitely generated $R$-modules, $V\stackrel{L}{\to}W$.
Its image, $L(V)$, is an $R$-submodule of $W$. Its cokernel, $coker(L)=\quotients{W}{L(V)}$, is an $R$-module as well. The annihilator-of-cokernel ideal is defined as the support of the cokernel module:
\begin{equation*}
\ca_L:=ann(coker(L)):=\{f\in R:\ f\quotients{W}{L(V)}=\{0\},\ i.e.\ fW\sseteq L(V)\}.
\end{equation*}
Recall the classical relation \cite[Proposition 20.7]{Eisenbud-book}:
 for $L\in\Mat$  with $m\le n$ there holds  $\ca_L\supseteq I_m(L)\supseteq (\ca_L)^m$.

(Here  $I_m(L)=I_{max}(L)$ is the ideal of the maximal, i.e. $m\times m$, minors.)
In particular, for $m=1$: $\ca_L=I_1(L)$, the ideal is generated by all the entries of $L$.

\

By definition $\ca_L W\sseteq L(V)$. In many cases one has the stronger property: $\ca_L W\sseteq J L(V)$, for some proper ideal $J\subsetneq R$.
\bex
Let $L\in Mat(m,n;J)$, $1<m\le n$ and suppose $\ca_L=I_{max}(L)$, e.g. this holds when $I_{max}(L)$ is radical.
Then $\ca_L W\sseteq J^{m-1} L(V)$.
\eex
The embedding $\ca_L W\sseteq\cm L(V)$ does not hold only in some degenerate cases. For example, let $L=\bpm f&0\\0&L_1\epm$, where $\det(L_1)=f$. Then $\ca_L=(f)$ and $(f)W\not\sseteq \cm\cdot L(V)$.

\section{The main results: criteria of solvability}
\subsection{The general statements}\label{Sec.Criteria.General}
Let $V\stackrel{L}{\rightarrow}W$ be a homorphism of (arbitrary) groups, where $V$ is filtered by normal subgroups as in \S\ref{Sec.Background.Filtered.Groups}.  Consider the equation $L(\uy)H(\uy)=L(\uv)$.
 See \S\ref{Sec.Background.Higher.Order.Terms} for the definitions.
\bthe\label{Thm.IFT.Nec.Suff.Criterion.non.Abelian.Groups}
1. If  the map $V\stackrel{H}{\to}W$ represents the ``higher order terms", i.e. $H(y)^{-1}H(yV_j)\sseteq L(V_{j+1})$ for any $y\in V_1$, $j\in\N$, then there exists a quasi-good \obo\ solution, $V_1\stackrel{y^{(n)}}{\to}V_1$.
 If moreover $L$ admits a right inverse,  a map $L(V)\stackrel{T}{\to}V$ satisfying $L\circ T=\one_{L(V)}$, $T(L(V_i))\sseteq V_i$, then there exists a good \obo\ solution.
\\2.  Suppose $V\stackrel{H}{\to}W$ is compatible with the filtration in the sense: $H(y\cdot V_j)\sseteq H(y)\cdot L(V_{N(j)})$, for some $N(j)$ satisfying $\liml_{j\to\infty}N(j)=\infty$.
If there exists a good \obo\ solution, $V_1\stackrel{y^{(n)}}{\to}V_1$, then $H$ represents the ``higher order terms", i.e. $H(yV_j)\sseteq H(y)\cdot L(V_{j+1})$ for any $y\in V_1$, $j\in\N$.
\\3. If $V$ is complete \wrt  $V_\bullet$ and $H$ represents the ``higher order terms" then there exists a quasi-good solution $V_1\stackrel{y}{\to}V_1$.
If moreover $L$ admits a right inverse, $L\circ T=\one_{L(V)}$, $T(L(V_i))\sseteq V_i$, then there exists a good  solution.
\ethe
\bpr
{\bf Part 1.} First we construct a quasi-good \obo\ solution $\uy^{(n)}$. The procedure is inductive with non-canonical choices. If $L$ is right-invertible then all the choices are canonical and the solution becomes good.

 Note that $H(V_{i})\sseteq L(V_{i+1})$, cf. \S\ref{Sec.Background.Higher.Order.Terms}.
Fix some $v\in V_i$, we construct inductively the sequence $y^{(n)}$ satisfying $y^{(n+1)}(y^{(n)})^{-1}\in V_{i+n}$ and  $L(y^{(n)})H(y^{(n)})L(v)^{-1}\in L(V_{i+n})$.

Choose $\uy^{(1)}=v$ and note that $L(\uy^{(1)})H(\uy^{(1)})L(v^{-1})=L(v)H(v)L(v)^{-1}\in H(V_i)\sseteq L(V_{i+1})$.

Suppose $\uy^{(1)},\dots,\uy^{(n)}$ have been constructed for some $n\ge1$. Present $\uy^{(n+1)}=z\uy^{(n)}$, so we should find the necessary $z\in V_{i+n}$.
 Note:
 \begin{multline}\label{Eq.proof.main.thm.formula}
L(\uy^{(n+1)})H(\uy^{(n+1)})L(\uv^{-1})=L(z)L(\uy^{(n)})H(\uy^{(n)})H(\uy^{(n)})^{-1}H(z\uy^{(n)})L(\uv)^{-1}
=\\=
\Big(L(z)\underbrace{L(\uy^{(n)})H(\uy^{(n)})L(\uv)^{-1}}_{w}\Big)L(\uv)H(\uy^{(n)})^{-1}H(z\uy^{(n)})L(\uv)^{-1}.
\end{multline}
By the induction assumption $w\in L(V_{i+n})$, thus we choose $z\in V_{i+n}$ satisfying $L(z)=w^{-1}$.
 (If $w$ is the identity element of $W$ then we choose $z={\bf \one}_V$.)
 Then equation \eqref{Eq.proof.main.thm.formula} reads:
\begin{equation*}
L(\uy^{(n+1)})H(\uy^{(n+1)})L(\uv^{-1})=L(\uv)\underbrace{H(\uy^{(n)})^{-1}H(z\uy^{(n)})}_{\in L(V_{i+n+1})}L(\uv^{-1})\in  L(V_{i+n+1}).
\end{equation*}
This completes the induction step. (Here we use the normality $V_j\lhd V$.)

By construction $y^{(n)}$ is a Cauchy sequence, as $y^{(n+1)}(y^{n})^{-1}\in V_{i+n}$. And if $v={\bf 1}_V$ then $y^{(n)}={\bf 1}_V$. Moreover, if  $v\in V_i$ then $y^{(n)}_v v^{-1}\in V_{i+1}$. Thus $y^{(n)}$ is a quasi-good \obo\ solution.

\

Suppose there exists a continuous right inverse, $L\circ T=\one_{L(V)}$, then in equation \eqref{Eq.proof.main.thm.formula} we choose \\$z_v=T(L(\uv)H(\uy^{(n)}_v)^{-1}L(\uy^{(n)}_v)^{-1})$. The proof of
$y^{(n)}_{v\De_j}(y^{(n)}_v\De_j)^{-1}\in V_{j+1}$ goes by induction on $n$. For $y^{(1)}_v=v$ the statement is trivial. Suppose this holds for $y^{(n)}_v$. Then
\begin{equation*}
y^{(n+1)}_{v\De_j}(y^{(n+1)}_v\De_j)^{-1}=(z^{(n+1)}_{v\De_j}y^{(n)}_{v\De_j})(z^{(n+1)}_{v}y^{(n)}_v\De_j)^{-1}.
\end{equation*}
Note that $y^{(n)}_{v\De_j}(y^{(n)}_v\De_j)^{-1}\in V_{j+1}$, thus
\begin{equation*}
(z^{(n+1)}_{v\De_j}y^{(n)}_{v\De_j})(z^{(n+1)}_{v}y^{(n)}_v\De_j)^{-1}\in
z^{(n+1)}_{v\De_j}V_{j+1}(z^{(n+1)}_{v})^{-1}.
\end{equation*}
 Now by the normality ($V_{j+1}\lhd V$):
\begin{equation*}
z^{(n+1)}_{v\De_j}V_{j+1}(z^{(n+1)}_{v})^{-1}=T\Big(L(\uv\De_j)H(\uy^{(n)}_{v\De_j})^{-1}L(\uy^{(n)}_{v\De_j})^{-1}\Big)
T\Big(L(\uv)H(\uy^{(n)}_v)^{-1}L(\uy^{(n)}_v)^{-1}\Big)^{-1}V_{j+1},
\end{equation*}
 which by normality ($V_1\lhd V$) equals
\begin{equation*}
T\Big(L(\De_j)\underbrace{H(\uy^{(n)}_{v\De_j})^{-1}}_{\in H(\uy^{(n)}_{v})^{-1}V_{j+1} }L(\uy^{(n)}_{v\De_j})^{-1}L(\uy^{(n)}_v)H(\uy^{(n)}_v)\Big)V_{j+1}.
\end{equation*}
Finally, by induction this space is $T\Big(L(\De_j)H(\uy^{(n)}_v)^{-1}L(\De_j)^{-1}H(\uy^{(n)}_v)\Big)V_{j+1}$.
And by the property of the filtration, $[V_1,V_j]\sseteq V_{j+1}$, the last expression is $V_{j+1}$

{\bf Part 2.} In Step 1 we prove that a good \obo\ solution is an almost surjective map, its image is dense. In Step 2 we use this auxiliary statement to bound $H^{-1}(y)H(y\De_j)$.

{\em Step 1.} We prove an auxiliary statement:
\begin{equation*}\ber
\text{\em if $V_1\stackrel{y}{\to}V_1$ is a good map,
 i.e. $y_v=v g(v)$ with $g(vV_j)g^{-1}(v)\sseteq V_{j+1}$,  then the image of $y$}
\\\text{\em  is dense in $V_1$, i.e. for any $v\in V_1$ there exists a sequence $v^{(n)}\in V_1$ such that $y_{v^{(n)}}v^{-1}\in V_{n+1}$.}
\eer\end{equation*}

Define $v^{(1)}=v$. Suppose $v^{(1)},\dots,v^{(i)}$ have been constructed, so that $y_{v^{(i)}}v^{-1}\in V_{i+1}$, i.e. $y_{v^{(i)}}v^{-1}=\De_{i+1}$. Define $v^{(i+1)}=\De_{i+1}^{-1} v^{(i)}$. Now the direct check:
\begin{equation*}
y_{v^{(i+1)}}=y_{\De_{i+1}^{-1}v^{(i)}}=\De_{i+1}^{-1}v^{(i)}\cdot g(\De_{i+1}^{-1}v^{(i)})\in \De_{i+1}^{-1}v^{(i)}\cdot g(v^{(i)})\cdot V_{i+2}=
\De_{i+1}^{-1}y_{v^{(i)}}\cdot V_{i+2}=v\cdot V_{i+2},
\end{equation*}
i.e. $v^{(i+1)}$ satisfies the needed condition.

\

{\em Step 2.} Fix some  good \obo\ solution $V_1\stackrel{y^{(n)}}{\to}V_1$, $L(y^{(n)}_v)H(y^{(n)}_v)\in L(v)L(V_{n+1})$. We should bound $H^{-1}(y)H(y\De_j)$, for any $y\in V_1$, $\De_j\in V_j$.
By Part 1 we can assume  $L(y)\in L(y^{(n)}_v V_{n+1})$ for some $v\in V_1$ and $n>j$. Moreover, we can choose $n$ so large that in addition: $L(y\De_j)\in L(y^{(n)}_{v\tilde{\De_j}} V_{n+1})$ for some $\tilde\De_j\in V_j$.
Now use $H(yV_j)\sseteq H(y)L(V_{N(j)})$ and choose $n$ so large that
 $H(y)\in H(y^{(n)}_v)L(V_{j+1})$ and  $H(y\De_j)\in H(y^{(n)}_{v\tilde{\De_j}})L(V_{j+1})$.
Therefore, for $n>j$:
\begin{multline}
H^{-1}(y)H(y\De_j)\in H^{-1}(y^{(n)}_{v})H(y^{(n)}_{v\tilde{\De_j}})L(V_{j+1})=
\Big(L(y^{(n)}_{v})^{-1}L(v)\Big)^{-1}\Big(L(y^{(n)}_{v\tilde{\De_j}})^{-1}L(v\tilde{\De_j})\Big)L(V_{j+1})=\\
=L\Big(v^{-1}y^{(n)}_{v}(y^{(n)}_{v\tilde{\De_j}})^{-1}v\tilde{\De_j}\Big)L(V_{j+1})=
L\Big(g(v)g(v\tilde{\De_j})^{-1}\Big)L(V_{j+1})\sseteq L(V_{j+1})
\end{multline}
In the second row we used the goodness of $y^{(n)}_v$.

\

{\bf Part 3.} If $V$ is complete \wrt $V_\bullet$ then $\capl_{i\ge1}V_i=\{{\bf 1}_V\}$. Given the Cauchy sequence $\uy^{(n)}$ from part 1, take the limit $\uy=\liml_{n\to\infty}\uy^{(n)}\in V$.
 Then one has
 \begin{equation*}
L(y)H(y)L(v)^{-1}=\liml_{n\to\infty}L(\uy^{(n)})H(\uy^{(n)})L(v)^{-1}\in\capl_n L(V_n)=\{{\bf 1}_W\}.\text{\epr}
 \end{equation*}
\beR
To emphasize, this theorem is almost an iff statement, thus the assumptions on $L,H$ are the ``weakest possible".
\eeR

\subsubsection{The case of abelian groups}
One often needs  results of such type for abelian groups, where one solves the equation $L(\uy-\uv)+H(\uy)=0$.  We state the corresponding criterion separately.
\bcor\label{Thm.IFT.Suff.Criterion.Abelian.Groups}
1. Given abelian groups $V$, $W$ and a homomorphism $V\stackrel{L}{\to}W$.
Suppose there exists a decreasing filtration $V_\bullet$ of $V$ satisfying:
 $\forall \uy\in V_1$, $\forall j\ge1$: $H(\uy+V_j)-H(\uy)\sseteq L(V_{j+1})$. Then for any $\uv\in V_1$ there exists a quasi-\wt\ \obo\ solution.
\\2. If $V$ is complete \wrt $V_\bullet$ then the conditions
 $\{H(\uy+V_j)-H(\uy)\sseteq L(V_{j+1})\}_{j\ge1}$ imply a quasi-\wt\ solution of the equation $L(\uy-\uv)+H(\uy)=0$.
\ecor

\beR
In the classical case of the equation $F(\bx,\uy)=0$ one asks that the map $V\stackrel{F'_\uy(\bx,0)}{\to}W$ is right invertible, i.e. surjective. Our criterion asks that
 $F'_\uy(\bx,0)(V_{j+1})$ contains the variation of the higher order terms, $H(y+V_j)-H(y)$, here
 $H(y)=F(\bx,\uy)-F(\bx,0)-F'_\uy(\bx,0)\uy$.
\eeR

\subsubsection{A special version for group-action equations}\label{Sec.Group.Action.IFT}
Given two maps of (not necessarily abelian) groups, $V\stackrel{L}{\underset{F}{\rightrightarrows}}W$. Suppose $W$ is filtered by the normal subgroups $\{W_j\}$ and $F(V)$ contains ${\bf 1}_W\in W$.
 (Here we do not assume that $L$ is a homomorphism.)
 Denote by $\overline{L(V)},\overline{F(V)}\sseteq W$ the closures \wrt the filtration $W_i$.  The following statement is almost tautological, yet highly useful in \S\ref{Sec.Group.Action.Equations.applications}.
\bel\label{Thm.Group.Action.IFT}
Suppose $F(V)\cdot\Big(L(V)\cap W_j\Big)\sseteq F(V)\cdot W_{j+1}$ for any $j\ge k$. If $W_k\sset \overline{L(V)}$ then $W_k\sseteq\overline{F(V)}$.
\eel
In the abelian case the condition reads: $F(V)+\Big(L(V)\cap W_j\Big)\sseteq F(V)+W_{j+1}$.
\bpr
Suppose $W_k\sseteq L(V)$, then $W_j\sseteq L(V)$ for any $j\ge k$. Thus $F(V)\cdot \Big(L(V)\cap W_k\Big)\sseteq F(V)\cdot W_{k+1}\sseteq\cdots$. As $F(V)$ contains ${\bf 1}_W\in W$ we get: $W_k\sseteq F(V)\cdot  W_N$ for any $N$.
Which is precisely $W_k\sseteq\overline{F(V)}$.

The general case. Let $N>k$, consider the quotient $W\stackrel{\pi_N}{\to}\quotients{W}{W_N}$. Denote the composition maps  $V\stackrel{L}{\underset{F}{\rightrightarrows}}W\stackrel{\pi_N}{\to}\quotients{W}{W_N}$ by $\pi_N L$, $\pi_N F$. Then $W_k\sset \overline{L(V)}$ implies $\pi_N(W_k)\sseteq \pi_N L(V)$ for any $N$. By the previous paragraph we get
 $\pi_N(W_k)\sseteq \pi_N F(V)$. Thus $W_k\sseteq F(V)\cdot W_N$ for any $N$. Which means $W_k\sseteq\overline{F(V)}$.
\epr

\subsection{Criteria for modules over the rings}
Theorem \ref{Thm.IFT.Nec.Suff.Criterion.non.Abelian.Groups} and corollary \ref{Thm.IFT.Suff.Criterion.Abelian.Groups} transform the solvability question into the search for the appropriate filtration $V_\bullet$.
 Not much can be said for a general (non)abelian group.  However our criterion simplifies
 for modules over a ring: it is enough to find just the first submodule $V_1\sset V$ and an ideal.

Let $R$ be a (commutative, associative) ring over a domain $\k$ of zero characteristic (e.g. $\k$ is a field).
Given two $R$-modules and a homomorphism, $L\in Hom_R(V,W)$. Suppose further that the term
 $H(\uy)$ admits a ``linear approximation  with the remainder in the form of Lagrange":
\beq\label{Eq.Taylor.Expansion.of.H}
H(\uy+\De)=H(\uy)+H_1(\uy)(\De)+H_2(\uy,\De)(\De,\De),
\eeq
here $H_1(\uy)(z)$ is linear in $z$ while $H_2(\uy,\De)(z,z)$ is quadratic in $z$.
\bex\label{Ex.Rings.with.Taylor.approximation}
Such an approximation holds e.g. for $R$ a subring of one of the quotients $\quotients{\k[[\bx]]}{I}$, $\quotients{C^\infty(\R^m,0)}{I}$.
\eex
\bcor\label{Thm.IFT.modules.over.rings.with.Taylor}
 Fix some ideal $J\sset R$ and a submodule $V_1\sset V$. Under the assumptions as above:
\\1. If the equation $L(y-v)+H(y)=0$ admits a good \obo\ solution for the filtration $\{V_i:=J^{i-1}V_1\}$ then $H(V_1)\sseteq J\cdot L(V_1)$.
\\2. If  $H(V_1)\sseteq J\cdot L(V_1)$ then $L(y-v)+H(y)=0$ admits a quasi-good \obo\ solution for the filtration $\{V_i:=J^{i-1}V_1\}$. (If $L$ is right invertible then there exists a good \obo\ solution.)
\ecor
\bpr
1. By theorem \ref{Thm.IFT.Nec.Suff.Criterion.non.Abelian.Groups} the existence of a \wt\ solution implies
 $H(\uy+\De_j)-H(\uy)\in J^j L(V_1)$ and hence $H(V_1)\sseteq J L(V_1)$.

2. For any $t\in\k$
 and $\De\in V_1$ we have $H(y+t\De)-H(y)\in JL(V_1)$. Thus
 $tH_1(\uy)(\De)+t^2H_2(\uy,t\De)(\De,\De)\in JL(V_1)$, for $t\in\k$. Then
$H_1(\uy)(\De)\in JL(V_1)$ and $H_2(\uy,t\De)(\De,\De)\in JL(V_1)$. Thus $H_1(\uy)(J^k\De)\sseteq J^{k+1}L(V_1)$ and $H_2(\uy,t\De)(J^k\De,J^k\De)\sseteq J^{k+1}L(V_1)$.
 Thus  $H(\uy+\De_j)-H(\uy)\in J^j L(V_1)$ is implied by $H(V_1)\sseteq J L(V_1)$. Now invoke corollary
 \ref{Thm.IFT.Suff.Criterion.Abelian.Groups} for the filtration $\{V_i:=J^{i-1}V_1\}$.
\epr

Corollary \ref{Thm.IFT.modules.over.rings.with.Taylor} reduces  the question (for modules over a ring) to the search for an appropriate submodule $V_1\sset V$.
 The simplest submodule is $V_1=JV$, for some ideal $J\sset R$.
\bcor
Suppose $H(y)$ has the linear approximation, as in equation \eqref{Eq.Taylor.Expansion.of.H}, and moreover $H(y)$ is of order $k\ge2$, i.e. $H(JV)\sseteq J^k W$. If $J H(V)\sseteq J\cdot L(V)$ then there exists a quasi-good \obo\ solution $J \cdot V\stackrel{y}{\to}J\cdot V$ \wrt the filtration $\{J^i V\}$.
\ecor
\noindent(proof: Apply corollary \ref{Thm.IFT.modules.over.rings.with.Taylor} for the filtration $V_i=J^i V$.)
\bex\label{Ex.obo.solution.for.Higher.Order.Terms}
Given the equation $L(y-v)+H(y)=0$, where $H$ is of order $k$, as above.

1. Consider the annihilator of cokernel ideal, $\ca_L:=ann(coker(L))$, cf. \S\ref{Sec.Background.Ann.Coker.and.J}.
\ls If $J^k\sseteq J\cdot\cm\cdot\ca_L$ then there exists a good (\obo) solution $J \cdot V\stackrel{y}{\to}J\cdot V$. \ls A bit weaker form: if
$J^k\sseteq J\cdot\ca_L$ then there exists a good (\obo) solution $\cm\cdot J \cdot V\stackrel{y}{\to}\cm\cdot J\cdot V$.

In the lowest order case, $k=2$, we get a sufficient condition for the \obo\ solvability: $J^2W\sseteq\cm JL(V)$.
This condition is weaker than Tougeron's and Fisher's conditions, so even this criterion is stronger.

2. Quite often $\cm\cdot L(V)\supseteq\ca_L\cdot W$, cf. \S\ref{Sec.Background.Ann.Coker.and.J}. Then we get a stronger statement: if
  $J^k\sseteq J\cdot\ca_L$ then there exists a good (\obo) solution $J \cdot V\stackrel{y}{\to}J\cdot V$.
\eex

\subsubsection{Ideals that satisfy $J^2\sseteq J\ca_L$}\label{Sec.Background.Ideals.Satisfying.J^2=..}
(These are important in view of  example \ref{Ex.obo.solution.for.Higher.Order.Terms}.)
Consider the set $\mathfrak{J}$ of all the ideals satisfying $J^2\sseteq J\ca_L$.
This is an inductive set, i.e. for
 any increasing sequence, $J_1\sseteq J_2\sseteq\cdots$ the union $\cupl_k J_k$
 is an ideal that satisfies $J^2\sseteq J\ca_L$. (If $f,g\in\cupl_k J_k$ then $f,g\in J_k$
  for some $k<\infty$, thus $fg,f+g\in J_k$.)
 Therefore in $\mathfrak{J}$ there exist(s) ideal(s) that is/are maximal by inclusion.

\bel Let $J\sset R$ be a maximal by inclusion ideal that satisfies $J^2\sseteq J\ca_L$.
\\1. $\ca_L\sseteq J$. If $J$  is finitely generated then $J\sseteq\overline{\ca_{L}}$.
(Here $\overline{\ca_{L}}$ is the integral closure.)
\\2. If $\ca_L$ is radical then $J=\ca_L$. If $R$ is integrally closed and $\ca_L$ is principal,
 generated by a non-zero divisor, then $J=\ca_L$.
\eel
\bpr
{\bf 1.} If $J^2\sseteq J\ca_L$ then obviously the inclusion is satisfied by the ideal $J+\ca_L$ as well. As $J$  is the largest with this property, $\ca_L\sseteq J$.

For the second part, note that $\ca_L$ is a reduction of $J$, see \cite[Definition 1.2.1]{Huneke-Swanson},
thus $J\sseteq\overline{\ca_{L}}$ by \cite[Corollary 1.2.5]{Huneke-Swanson}.

{\bf 2.} If $J^2\sseteq J\ca_L$ then in particular $J^2\sset\ca_L$. Then, $\ca_L$ being radical, we get $J\sseteq\ca_L$. Together with part 1 we get: $J=\ca_L$.

The second part follows from \cite[Proposition 1.5.2]{Huneke-Swanson}: in our case $\overline{\ca_{L}}=\ca_L$.
\epr
\bex In many cases $\ca_L\subsetneq J\subsetneq\overline{\ca_{L}}$ and a maximal by inclusion ideal $J$ is non-unique.  For example, let $R=\k[[x,y,z]]$ and $L=(x^p,y^p,z^p)\in Mat(1,3;R)$. Then $\ca_L=(x^p,y^p,z^p)$ while
$\overline{\ca_{L}}=\cm^p$.
 Define: $J_z=((x,y)^p,z^p)$, $J_y=((x,z)^p,y^p)$, $J_x=((y,z)^p,x^p)$. By the direct check,
  each of them satisfies $J^2=J\ca_L$. But there is no bigger ideal $J$ that contains say $J_x+J_y$
  and satisfies $J^2=J\ca_L$. Indeed, suppose $y^{p-i}z^{i}\in J$ and $x^{p-j}z^{j}\in J$, for some $i,j$
  satisfying $i+j<p$. Then $J\ca_L=J^2\ni x^{p-j}y^{p-i}z^{i+j}$, in particular $x^{p-j}y^{p-i}z^{i+j}\in\ca_L$ for  $i+j<p$, contradicting the definition of $\ca_L$.
Thus, in this case there are at least 3 distinct maximal by inclusion ideals.
\eex

\subsection{(Non-)Uniqueness}\label{Sec.Uniqueness}
The classical Implicit Function Theorem ensures the uniqueness of solution, provided $F'_\by(0,0)$ is injective. In our case the injectivity ensures that the solution is ``eventually unique" in the following sense.
\bprop
Given two \obo-solutions $y^{(n)}_1$, $y^{(n)}_2$ of the equation $L(y)H(y)=L(v)$.
Suppose $y^{(1)}_1$, $y^{(1)}_2\in V_1$ and $L$ is injective. Then for any $n$: $y^{(n)}_1(y^{(n)}_2)^{-1}\in V_n$.
\eprop
\bpr
By the assumption $y^{(1)}_1(y^{(1)}_2)^{-1}\in V_1$. Suppose the statement holds for $j=1,\dots,(n-1)$.
As both $y^{(n)}_i$ are Cauchy sequences we get  $y^{(n)}_1(y^{(n)}_2)^{-1}\in V_{n-1}$. We prove that in fact
$y^{(n)}_1(y^{(n)}_2)^{-1}\in V_{n}$.

As each $y^{(n)}_i$ is an \obo-solution we have $L(y^{(n)}_i)H(y^{(n)}_i)\in L(v)V_n$. Thus
\begin{equation*}
\Big(L(y^{(n)}_2)H(y^{(n)}_2)\Big)^{-1}\Big(L(y^{(n)}_1)H(y^{(n)}_1)\Big)\in V_n.
\end{equation*}
By the normality,  $V_n\lhd V$, we get:
\begin{equation*}
\Big(L(y^{(n)}_2)\Big)^{-1}L(y^{(n)}_1)\in H(y^{(n)}_2)\Big(H(y^{(n)}_1)\Big)^{-1}V_n.
\end{equation*}
  Now use $y^{(n)}_1(y^{(n)}_2)^{-1}\in V_{n-1}$ and the property of higher order terms for $H$ to get:
 $H(y^{(n)}_2)\Big(H(y^{(n)}_1)\Big)^{-1}\in L(V_n)$. Therefore $L\big((y^{(n)}_2)^{-1}y^{(n)}_1\big)\in L(V_n)$ and the statement follows by the injectivity of $L$.
\epr
\beR
The assumption $y^{(1)}_1$, $y^{(1)}_2\in V_1$ is important. One might seek for a condition in terms of $\uv$ and $L$ only, then it is natural to ask that $\uv$ belongs to a small enough subgroup of $V$. For example, in the case of modules, $v\in JV$, for some small enough ideal $J\sset R$.
This does not suffice as one sees already in the example of one equation in one variable: $(y-x^a)(y+x^b)=0$.
 Suppose $a<b$, then $\ca_L=(x^a)$, while $v\in (x^{a+b})$. By taking $b\gg a$ the ideal $(x)^{a+b}$ can be made arbitrarily small as compared to $\ca_L$. Yet, there is no uniqueness.
\eeR
\beR
If $L$ is non-injective then there can be no uniqueness. Even the images $L(y^{(n)})$  of an \obo-solution are not ``eventually unique". As the simplest example consider the equation $y^2_1+y^2_2-y_1+v=0$, where $v\in(x)\sset\k[[x]]$. We have a family of solutions $y_1=\frac{2(v+y^2_2)}{1+\sqrt{1-4(v+y^2_2)}}$, here $y_2$ is a parameter. By taking $y_2\in(v^j)$ these solutions can be made arbitrarily close one to the other (in particular they all lie in $V_1$), yet $L(y_1,y_2)$ is different for different $y_2$.
\eeR

\subsection{A criterion for exact solutions.}\label{Sec.Exact.Solutions.Over.Rings.with.IFT}
The criteria of \S\ref{Sec.Criteria.General} provide \obo\ solutions, alternatively: solutions in the completion of $V$ by $V_\bullet$, i.e. the formal solutions.
Recall the Artin approximation property: {\em if a finite system of polynomial equations over $R$ has a solution over $\hat{R}$ then it has a solution over $R$}, \cite{Artin68}.
Many rings have this approximation property, for example excellent Henselian rings (in particular complete rings, analytic rings), cf. \cite{Hauser-Rond}.

In our case we have more general rings and more general class of equations. Thus we give a criterion for exact (and not just \obo) solution.

 Fix some proper ideal $J\ssetneq R$. The pair $(R,J)$ is said to satisfy the (classical) implicit function theorem, denote this by $cIFT_J$, if for any {\em surjective} morphism of free $R$-modules of finite ranks, $V\stackrel{L}{\to}W$, any $v\in J V$ and any ``higher order term", $V\stackrel{H}{\to}W$, the equation $L(y-v)+H(y)=0$ has a good solution.
  Note that if $R$ satisfies $cIFT_J$ then for any ideal $J_1\sseteq J$ the ring satisfies $cIFT_{J_1}$ as well.
\bex\label{Ex.Rings.Satisfying.classical.IFTm}
Let $(R,\cm)$ be any local Henselian ring over a field $\k$. For example, the ring of formal power series $R=\quotients{\k[[x_1,\dots,x_m]]}{I}$, the ring of analytical power series
 $R=\quotients{\k\{x_1,\dots,x_m\}}{I}$ (for $\k$-normed),  the ring of smooth functions $R=\quotients{C^\infty(\R^m,0)}{I}$ or the ring of $p$-times differentiable functions  $R=\quotients{C^p(\R^m,0)}{I}$.
Then $(R,\cm)$ satisfies the $cIFT_\cm$.

The rings $\k[\bx]$, $\k[\bx]_{(\bx)}$ do not satisfy $cIFT_\cm$, e.g. the equation $y^2+y=x^2$ is not solvable over these rings.
\eex
We say that $(R,J)$ satisfies the {\em implicit function theorem with unit linear part}, denote this by $IFT_{J,\one}$, if the system of equations $\uy-\uv+H(\uy)=0$ has a good solution $JV\stackrel{y}{\to}JV$,
for any higher order terms $H$.

This system is a particular case of the classical implicit function equations. Therefore the Henselian rings (over a field) of example \ref{Ex.Rings.Satisfying.classical.IFTm} satisfy $IFT_{\cm,\one}$.
Note that the condition $IFT_{J,\one}$ is  weaker than $IFT_J$.
 For example, $IFT_{J,\one}$ is satisfied by $\quotients{\Z[[\bx]]}{I}$, $\quotients{\Z\{\bx\}}{I}$, for $J=(\bx)$.
 More generally, one can take $\k$ any ring and $R$ a Henselian algebra over $\k$.

\bprop\label{Thm.I.F.T.Suff.Crit.Modules}
Given a finitely generated $R$-module $V$ and two maps $V\stackrel{L}{\underset{H}{\rightrightarrows}}W$. Suppose
 $L\in Hom_R(V,W)$, while $H$ satisfies $H(\suml_i y_i\xi_i)=\suml_i h_i(\{y_j\})L(\xi_i)$, here $\{\xi_i\}$ are some  generators of $V$, while $h_i(\{y_j\})$ are of order$\ge2$. Suppose $IFT_{J,\one}$ holds for an ideal $J\ssetneq R$. Then the equation $L(y-v)+H(y)=0$ has a solution $JV\stackrel{y}{\to}JV$.
\eprop
\noindent Note that here $R$ is not necessarily over a field, e.g $R$ can be $\Z[[\bx]]$ or $\Z\{\bx\}$. Being of order$\ge2$ means that $h_i(J)\sseteq J^2h(R)$ for any ideal $J\sseteq R$.
\bpr
Expand $v=\suml_i v_i\xi_i$, $y=\suml_i y_i\xi_i$, then the equation reads: $\suml_i(y_i-v_i+h_i(y))L(\xi_i)=0$.
Thus it is enough to solve the finite system of equations $\{y_i-v_i+h_i(y)=0\}$. As $IFT_{J,1}$ holds in our situation we get the solution.
\epr

\bcor\label{Thm.IFT.J^2=JA}
Suppose a local ring $(R,\cm)$ satisfies $IFT_{\cm,\one}$. Consider the equation $L(\uy-\uv)+H(\uy)=0$, where $ord(H)\ge2$.
\\1. If $J^2W\sseteq\cm J L(V)$ then for any $v\in J V$ there exists a  solution.
\\2. If $J^2\sseteq J\ca_L$ then for any  $v\in \cm J V$ there exists a solution.
\\3. If $J^2\sseteq J\ca_L$ and $\ca_L W\sseteq\cm L(V)$ then for any  $v\in  J V$ there exists a  solution.
\ecor
\bex
Let $(R,\cm)$ be a local Henselian ring over a field. Take $J=\ca_L$, then the corollary implies
 Tougeron's and Fisher's theorems. As mentioned in the introduction, if one takes $J$ the maximal possible that
  satisfies $J^2=J \ca_L$ then one gets the strengthening of  Tougeron's and Fisher's theorems.

But the corollary is useful for more general rings, e.g. if in equation \eqref{Eq.hard.equation}
the term $p(x_1,x_2)$ has integral coefficients then we get a solution over $\Z[[x_1,x_2]]$.
\eex

\section{Examples, remarks and applications}\label{Sec.Remarks}

\subsection{Comparison to Fisher's and Tougeron's theorems.}\label{Sec.Comparison.to.Fisher.Tougeron}
The condition $J^2\sseteq J\cdot\ca_{F'_{\uy}(\bx,0)}$ (cf. corollary \ref{Thm.IFT.J^2=JA}) is a weakening of the condition
  $J\sseteq \ca_{F'_{\uy}(\bx,0)}$.
\bex
Let $R=\k[[x_1,x_2]]$, where $\k$ is some base ring, take $\cm=(x_1,x_2)$. (If $\k$ is a field then $\cm$ is the  maximal ideal.)
Consider the equation $H(\by,\bx)+y_1x^k_1+y_2x^k_2+p(\bx)=0$, compare this to equation \eqref{Eq.hard.equation}.
 Here $H(\by,\bx)$ represents the higher order terms, it is at least quadratic in $y_1,y_2$, while  $p(\bx)\in R$.
 In this case: $L=F'_{\by}(\bx,0)=(x^k_1,x^k_2)\in Mat(1,2;R)$
 and $I_{max}(L)=\ca_{L}=(x^k_1,x^k_2)\sset R$.
  Thus $(\ca_L)^2=(x^{2k}_1,x^{k}_1x^{k}_2,x^{2k}_2)$.
 Thus to apply Tougeron's and Fisher's theorems we have to assume: $p(\bx)\in \cm(x^{2k}_1,x^{k}_1x^{k}_2,x^{2k}_2)$.
\\On the other hand,  by direct check, the ideal $\cm^k=(x_1,x_2)^k$ satisfies:
$(\cm^k)^2=\cm^k\cdot(x^k_1,x^k_2)=\cm^k\cdot\ca_{L}$. Therefore corollary \ref{Thm.IFT.J^2=JA} gives:
if   $p(\bx)\in \cm^{2k+1}$ then the equation has a solution.

For $\k$ an algebraically closed field we get a better criterion: if $p(\bx)\in \cm^{2k}$ then the equation has a solution.

\

Note that to write down an explicit solution is not a trivial task even in the particular case of equation \eqref{Eq.hard.equation}.

Further, if $\k$ is not a field then we get the solvability of a ``Diophantine type" equation. For example, for $\k=\Z$ and $H(y,\bx),p(\bx)$ defined over $\Z$, we get the criterion of solvability over $\Z[[\bx]]$. Note that even for the equation $y^n+yx^k+x^N=0$ the solvability over $\Z[[\bx]]$ is not totally obvious.
\eex

\

Therefore, even in the case of just one equation,  the condition $J^2=J\cdot\ca_L$ strengthens the versions of Tougeron and Fisher.

\subsection{Comparison of the condition $J^2=J\cdot\ca_{F'_{\uy}(\bx,0))}$ to $H(V_1)\sseteq J L(V_1)$}\label{Sec.Comparison.of.ideal.conditions.to.filtration} (cf. corollary \ref{Thm.IFT.modules.over.rings.with.Taylor})
It is simpler to check the ideals, $J^2=J\cdot\ca_{F'_{\uy}(\bx,0))}$, than to look for a submodule satisfying the needed  property. But the ``ideal-type" criterion is in general weaker than the criterion via $V_1$.

\bex
Consider the system $\Big\{\ber y^2_1+y_1x_1=x^n_1\\y^2_2+y_2x_2=x^m_2\eer\Big\}$ over $R=\k[[x_1,x_2]]$.
In this case the annihilator of cokernel ideal is principal,
$\ca_{F'_{\uy}(\bx,0))}=(x_1x_2)$, thus $J^2=J\cdot \ca_{F'_{\uy}(\bx,0)}$
implies $J=\ca_{F'_{\uy}(\bx,0)}$, see \S\ref{Sec.Background.Ideals.Satisfying.J^2=..}.
And $(x_1x_2)$ does not contain $x^n_1$, $x^m_2$, regardless of how big are $n$ and $m$.

Of course, the general criterion of corollary \ref{Thm.IFT.modules.over.rings.with.Taylor} suffices here.
 (One starts from $V_1=\bpm x_1 R\\x_2 R\epm$ and $J=(x_1,x_2)$.)
\eex

This is a good place to see in a nutshell why no weakening of $J^2=J\cdot\ca_{F'_{\by}(\bx,0)}$ in the form
 of some condition on ideals is possible.
\bex
 Consider the related system with a modified quadratic part:
$\Big\{\ber y^2_1+y_1x_1=x^n_1\\y^2_1+y_2x_2=x^m_2\eer\Big\}$.
While the previous system has obvious solutions for $n,m\ge2$, this system has no solutions in $R$.
 Indeed, from the second equation it follows that $y_1$ is divisible by $x_2$. Then the left hand side
  of the first equation must be divisible as well, contradicting the non-divisibility of the right hand side.
\eex
\bex
As a further illustration we consider the system $\Big\{\ber y^2_2+y_1a_1=b_1\\y^2_1+y_2a_2=b_2\eer\Big\}$, where
  $a_i,b_i\in \cm\sset R$, here $R$ is a regular local Henselian ring. Suppose $gcd(a_1,a_2)=1$, i.e. $(a_1)\cap(a_2)=(a_1a_2)$.
   Then $\ca_L=(a_1a_2)$ is a principal ideal and thus $J^2=J\ca_L$ implies $J=\ca_L$. Thus the approach via $J^2=J\cdot\ca_L$ gives:
\beq\label{Eq.inside.Example}
\text{   if $b_1\in \cm(a^2_1a_2)$, $b_2\in\cm(a_1a^2_2)$ then the system has a solution.}
\eeq
We check the approach via filtration. To invoke the corollary \ref{Thm.IFT.modules.over.rings.with.Taylor} we need $V_1\sset R^{\oplus 2}$ to satisfy: if $\bpm v_1\\v_2\epm\in V_1$ then $\bpm v^2_2\\v^2_1\epm\in\bpm a_1v_1\\a_2v_2\epm$ for any $\bpm v_1\\v_2\epm\in V_1$.
 This gives: $V_1\subset(a_1a_2)R^{\oplus2}$.
  Put $V_1=(a_1a_2)R^{\oplus2}$, this ensures $H(V_1)\sset \cm L(V_1)$.
  Note that $L$ has the obvious continuous right-inverse, $\bpm a_1 h_1\\a_2h_2\epm\stackrel{T}{\to}\bpm h_1\\h_2\epm$. Thus for $b_1\in (a^2_1a_2)$, $b_2\in (a_1a^2_2)$ the equation has a good \obo\ solution.
 This later condition is slightly weaker than that of equation \eqref{Eq.inside.Example}.
\eex

\beR 
Suppose the system of equations splits. Then it is natural to choose the split submodule: $V_1=V_{1,1}\oplus V_{1,2}\sset V$. (Note that the converse does not hold: decomposability of $V_1$ does not imply that the system splits in any sense. For example, all the modules of the type $V_1=JV$ are decomposable if $V$ is free of $rank>1$.)
The following questions are important:
\ls Suppose $L$ is block-diagonal. What are the conditions on $H$ so that we can choose $V_1=V_{1,1}\oplus V_{1,2}$?
\ls Formulate some similar statements for $L$ upper-block-triangular vs $V_1$ an appropriate extension of modules.
\eeR

\subsection{Equations whose solutions are not good}\label{Sec.Equations.with.no.good.solutions}
Often the ``simple" and ``most natural" solutions are not \wt\ (not even quasi-good) in our sense, moreover the (quasi-)good solutions do not exist at all.
\bex
Consider the equation $y^2=p(x)$ over $R=\k[[x]]$. Here $L=0$, while $H(y)\neq0$. Thus corollary  \ref{Thm.IFT.modules.over.rings.with.Taylor} claims that there are no good solutions. Explicitly: there does not exist a submodule $\{0\}\neq V_1\sset R$ such that for any $p(x)\in V_1$ there exists a solution $V_1\stackrel{y}{\to}V_1$, good in the sense of \S\ref{Sec.Background.Good.Solution}.
This can be seen directly, if $V_1\neq\{0\}$ then $x^{2N+1}\in V_1$ for $N\gg1$, and $y^2=x^{2N+1}$ has no solutions in $R$.

Of course,  by the direct check, for any $p(x)$ of even order there are solutions. But these solutions are not good.
\eex
\bex
Consider the equation $(y-g_1(\bx))(y-g_2(\bx))=0$, over $R=\k[[\bx]]$, $\bx=x_1,\dots,x_n$, $n>1$. Assume that $g_1(\bx)$, $g_2(\bx)$ are generic enough, in particular $g_1(\bx)\cdot g_2(\bx)\not\in(g_1(\bx)+g_2(\bx))=L(V)$. Then the equation cannot be presented in the form $L(y-v)+H(y)=0$, so it has no quasi-good solutions.
(Even its linear part is non-solvable, though the equation has two obvious solutions.) This happens because an arbitrarily small deformation of the free term, $g_1(\bx)g_2(\bx)$, will bring an equation with no solutions in $R=\k[[\bx]]$. (In the case $g_1(\bx),g_2(\bx)\in C^\infty(\R^p,0)$ even a deformation by a flat function will lead to an equation with no solutions.)
\eex

\subsection{An application: bifurcations of polynomial roots}\label{Sec.Applications.Bifurcations} Fix a polynomial $p(y)=\suml^d_{i=0}a_i y^i$. Suppose for a fixed tuple of the coefficients, $(a_0,\dots,a_d)$, the polynomial has only simple roots (of multiplicity one). Then under small deformations of the coefficients the roots deform smoothly.

The multiple roots cause bifurcations in general. Our results provide a sufficient condition that a particular root deforms (smoothly/analytically/\dots) under the change of parameters. More precisely, starting from the initial ring $R$ consider an extension $S$ of $R$ by one local  variable, e.g. $S=R[[t]]$ or $S=R\{t\}$ etc.
 Present the family of the equations in the form $a_0(t)+a_1(t)y+H(y,t)=0$.
We say that a root $y_0$ of the  initial equation {\em deforms} ( smoothly/analytically/\dots) if there exists a root $y(t)\in S$ satisfying $y(0)=y_0$.

To formulate the criterion we shift the variables $y\to y+y_0$, so that the (new) root of the initial  equation is $y=0$.
\bcor
1. (Tougeron) If $a_0(t)\in \big(ta_1^2(t)\big)$ then the root $y=0$ of the initial equation deforms with $t$.
\\2. (B.-K.) If $a_0(t)\in(t a_1(t))$ and $a_i(t)(a_0(t))^{i-1}\in t (a_1(t))^i$ for any $i\ge2$ then the root $y=0$ of the initial equation deforms with $t$.
\ecor
\noindent(Note that if $a_0(t)\in \big(ta_1^2(t)\big)$ then all the assumptions of part two are satisfied.)
\bpr
Here put $v=\frac{a_0}{a_1}$ and $V_j=t^{j-1}\cdot(v)$.
\epr
\bex
If all the eigenvalues of a matrix are distinct, then they deform differentiably under the small deformations of the entries. In the case of multiple eigenvalues the corollary above
 ensures  that at least one of the potentially bifurcating eigenvalues deforms differentiably. Explicitly, expand the determinant:
\begin{equation*}
 det(A_t-y\one)=\suml_{i\ge2}(-y)^itrace(\overset{n-i}{\wedge}A_t)=det(A_t)-trace(A^\vee_t)y+\cdots+(-y)^n.
\end{equation*}
(Here $\overset{n-i}{\wedge}A_t$ is the associated skew-power of $A_t$.)
Suppose the multiple eigenvalue is zero, so $\det(A_{t=0})=0$. Then
\ls(Tougeron's part) If $det(A_t)\in\Big(t\big(trace(A^\vee_t)\big)^2\Big)$ then the eigenvalue deforms smoothly.
\ls(B-K's part) If $trace(\overset{n-i}{\wedge}A_t)\cdot(\frac{det(A_t)}{trace(A^\vee_t)})^i\in\Big(t\cdot det(A_t)\Big)$ then the eigenvalue deforms smoothly.
\eex

\subsection{A possible application: smooth curve-germs on singular spaces}\label{Sec.Applications.Smooth.Germs}
Let $(X,0)\sset(\k^n,0)$ be a germ (algebraic/analytic/formal) of a singular space. The smooth curve-germs lying on $(X,0)$ is an important subject, often used in the theory of arc spaces, \cite{Denef-Loeser}. The first question is whether $(X,0)$ admits at least one smooth curve-germ, \cite{Gonzalez-Sprinberg.Lejeune-Jalabert}, \cite{J.O.}, \cite{J.O.D.S.brieskorn-pham}.

From the IFT point of view this question reads (for simplicity we work over $\k[[x_1,\dots,x_n]]$):
\[\ber
\text{\em Can a given system of equations be augmented by another system, so that}
\\\text{\em  the total system, $\{F(\bx,\by)=0=G(\bx,\by)\}$, has one-dimensional power series solutions?}
\\\text{\em For example, $x_2(x_1),\dots,x_n(x_1)\in\k[[x_1]]$, $F(x_1,x_2(x_1),\dots,x_n(x_1))\equiv0$}
\eer\]
The strong IFT seems to lead to some results on the existence/properties of families of such curves.

\section{An approximation theorem of Tougeron-Artin type}\label{Sec.Approximation.Theorem}
There are several approximation theorems guaranteeing analytic/$C^\infty$ solutions,
 provided a formal solution exists.
Given the germ of an analytic map at the origin, $F:(\R^m,0)\times(\R^n,0)\rightarrow(\R^p,0)$,
consider the Implicit Function Equation
\beq\label{Eq.Main.Equation}
F(\bx,\by)=0
\eeq
here $\bx$ is the multi-variable, while $\by$ is an unknown map, $(\R^m,0)\stackrel{\by}{\to}(\R^n,0)$.

A {\em formal solution} of this equation is a formal series $\hy(\bx)\in\R[[x]]^{\oplus n}$ satisfying:
$\hF(\bx,\hy(\bx))\equiv0$,
where $\hF$ is the (formal) Taylor expansion at zero of the map $F$. In general this solution
does not converge off the origin. Two classical results relate it to the ``ordinary" solution.
\bthe
 Let $\hy(\bx)$ be a formal solution of the analytic equation $F(\bx,\uy)=0$.
\\1. \cite{Artin68} For every $r\in\N$ there exists an
analytic solution whose $r$'th jet coincides with the $r$'th jet of $\hy(\bx)$.
\\2. \cite{Tougeron1976} There exists a $C^\infty$-solution $\uy(\bx)$, whose Taylor series at the origin
is precisely $\hy(\bx)$ and such that for any $r\in\N$ there exists an analytic solution which
 is $r$-homotopic to $\by(\bx)$.
\ethe
(Recall that two solutions, $\by_0(\bx)$, $\by_1(\bx)$, are $r$-homotopic if there exists a $C^\infty$
 family of solutions, $\by(\bx,t)$,  such that $\by_0(\bx)=\by(\bx,0)$, $\by_1(\bx)=\by(\bx,1)$,
 and $\by(\bx,t)-\by_0(\bx)$ is r-flat at the origin.)

What if the equation $F(\bx,\uy)=0$ is not analytic but only of $C^\infty$-type? Does the existence of a
formal solution, for the completion $\hF(\bx,\uy)=0$, imply the
 existence of a $C^\infty$ solution? The naive generalization of Artin's/Tougeron's theorems does not hold.
\bex
Let $\tau$ be a function flat at the origin, e.g. $\tau=\Big\{\ber e^{-\frac{1}{x^2}},\   x\neq0 \\0,\  x=0.\eer$.
 Consider the equation $ {\tau}^2(x)y(x) = {\tau}(x) $. The completion of this equation is the identity, $0\equiv0$,
 thus every formal series $\hy\in\R[[x]]$ is a formal solution of $\hF(\bx,\uy)=0$.
 However, the equation has no local smooth solutions (not even continuous ones).
\\In this example the coefficient of $y(x)$, i.e. the function $\tau^2$, is flat at zero.
 In other words, the ideal $\ca_{F'_\uy(\bx,\by_0)}$ is too small and $\ca_{F'_\uy(\bx,\by_0)}\cdot\cm^\infty\neq\cm^\infty$.
\eex
The following statement supplements our previous results, and extends Tougeron's theorem to $C^\infty$-equations.
Let $R=C^\infty(\R^m,0)$, with the maximal ideal $\cm\sset R$.
 Suppose the equation $F(\bx,\by)=0$ has a formal solution, $\hat{\by}_0$. By Borel's lemma, \cite{Rudin-book}, we can choose a $C^\infty$-map $\by_0$ whose completion is $\hat{\by}_0$, thus $F(\bx,\by_0)$ is a vector of flat functions.
\bthe\label{Thm.Approximation.Theorem.for.smooth.funcs}
Suppose the equation $F(\bx,\by)=0$ has a formal solution and $\det\Big(F'_\by(\bx,\by_0)(F'_\by(\bx,\by_0))^T\Big)\cm^\infty=\cm^\infty$.
 Then there exists a $C^\infty$-map, $(\R^m,0)\stackrel{\by}{\to}(\R^n,0)$,
  whose Taylor series at the origin
is precisely $\hy_0$ and   such that $F(\bx,\by(\bx))\equiv0$.
\ethe
\bpr
We seek the solution in the form $\by=\by_0 +\bz$, where the map $\bz$ is flat.
Expand $F(\bx,\by_0+\bz)$ into the Taylor series with remainder:
\begin{equation*}
F(\bx,\by_0+\bz) = F(\bx,\by_0) + F_\by'(\bx,\by_0)\bz+\left( \int\limits_0^1(1-t)\frac
{\partial^2F(\bx,\by_0+t\bz)}{\partial \uy^2}dt\right)(\bz)^2.
\end{equation*}
Then the equation  takes the form
\beq\label{Eq.Current}
F'_\by(\bx,\by_0)\bz + G(\bx,\bz) =-F(\bx,\by_0)
\eeq
where the map $F(\bx,\by_0)$ is flat. Note that the summand $G$ satisfies the condition
$G(\bx,\la \bz) = \la^2hH(\bx,\bz,\la)$
with a $ C^\infty$-map $H$ such that $ H_z'(\bx,0,\la) = 0$. We look for the solution of equation \eqref{Eq.Current}
 in the form
 \begin{equation*}
 z = d(\bx)\Big(F_\uy'(\bx,y_0)\Big)^T\Big(F_\uy'(\bx,y_0)(F_\uy'(\bx,y_0))^T\Big)^\vee u,\quad
 \end{equation*}
where $d(\bx):=\det\Big((F_\uy'(\bx,\by_0)(F_\uy'(\bx,\by_0))^T\Big)$ and $A^\vee $ denotes the adjugate matrix.

Then we arrive at the equation $d^2(\bx)u + d^{2}(\bx)\tilde{G}(\bx,u) =-F(\bx,\by_0)$
with the $C^\infty$-map $\tilde{G}$ satisfying $\tilde{G}_u'(\bx,0) = 0 $. Dividing by $ d^2(\bx) $, we obtain the equation
\begin{equation*}
u +\tilde{G}(\bx,u) = \tau(\bx),
\end{equation*}
where the map $\tau$ is flat.  By the classical Implicit Function Theorem, the latter
equation has a local flat $ C^\infty$-solution. Hence, the map
$ z $ satisfies the equation \eqref{Eq.Current}, and $\by =\by_0+\bz$ is the solution we need.
\epr
\beR
1.
The assumption of the theorem can be stated as:
\begin{center} {\em every function flat at the origin is divisible by
$\det\Big(F'_\uy(\bx,\by_0)(F'_\uy(\bx,\by_0))^T\Big)$.}
\end{center}
In particular this implies: $F'_\uy(\bx,\by_0)$ is non-degenerate in some punctured neighborhood of the origin $0\in (\R^m,0)$.

Note that $\by_0$ is defined up to a flat function, but the assumption does not depend on this  choice.

2.
Recall that a function $g(\bx)$ is said to satisfy the Lojasiewicz condition (at the origin) if there exist constants $C>0$ and $\de>0$ such that for any point $\bx\in(\R^m,0)$: $|g(\bx)|\ge C dist(\bx,0)^\de$.
 As is proved e.g. in \cite[\S V.4]{Tougeron-book}: $g(\bx)$ satisfies the Lojasiewicz condition at the origin iff
  $g(\bx)\cm^\infty=\cm^\infty$. Thus the assumption  of the theorem can be stated in the form:
\beq\label{Eq.Lojasiewicz.Condition}
\det\Big(F'_\uy(\bx,\by_0)F'_\uy(\bx,\by_0)^T\Big)\ge C dist(\bx,0)^\de,\ \text{for some $C,\de>0$}.
\eeq

3. A similar statement can be proved for $C^k(\R^p,0)$ functions, but then the solution is in general only in the $C^{k-2-\de}$ class.
\eeR

\section{Openness of group orbits and applications to the finite determinacy}\label{Sec.Group.Action.Equations.applications}

Given a module $W$ over some base ring $\k$ (we assume $\k\supseteq\Q$) with a decreasing filtration $\{W_i\}$.
 Consider the group of all the $\k$-linear invertible maps that preserve the filtration,
 $GL(W_\bullet):=\{g\in GL(W),\ g(W_j)=W_j\}$.
 Fix some subgroup $G\sseteq GL(W_\bullet)$ and let $G^0\sseteq G$ be the unipotent subgroup,  \S\ref{Sec.Group.Action.Definition.G0}. Fix some element $w\in W$, consider the germ of its $G^0$-orbit, $(G^0w,w)$ and the tangent space to this germ, $T_{(G^0w,w)}$, \S\ref{Sec.Group.Action.Logarithm.Exponent.Tangent.Space} (Note that the existence of $T_{(G^0w,w)}$ places some restrictions on $G$, see equation \eqref{Eq.Assumption.Group.Has.Logarithm}.)
\bthe\label{Thm.Determinacy}
If $W_k\sseteq\overline{T_{(G^0w,w)}}$ then $w+W_k\subseteq  \overline{G^0w}$.
\ethe
Here $\overline{(\cdots)}$ denotes the closure \wrt the filtration $W_\bullet$. Thus the statement is of the \obo-type. In particular, in the proof we can assume that $W$ is $W_\bullet$-complete.

The proof  is given in \S\ref{Sec.Group.Action.Proof.of.Theorem}, after some preparations in \S\ref{Sec.Group.Action.Preparations}. Some immediate applications to the finite determinacy are given in \S\ref{Sec.Group.Action.Determinacy.Appplications}.

\subsection{Preparations}\label{Sec.Group.Action.Preparations}
\subsubsection{The unipotent subgroup $G^0$}\label{Sec.Group.Action.Definition.G0} Consider the system of projecting maps $W\stackrel{\pi_j}{\to}\quotients{W}{W_j}$. They induce projections on the group $GL(W)\stackrel{\pi_j}{\to}GL(\quotients{W}{W_j})$ and accordingly the restrictions $G\stackrel{\pi_j}{\to}\pi_j(G)\sseteq GL(\quotients{W}{W_j})$. (We use the same letter $\pi_j$, this causes no confusion.) We define the ``unipotent" part of the group, $G^0:=\{g\in G,\ \forall j>0:\ \pi_jg|_{\quotients{W_{j-1}}{W_j}}=Id|_{\quotients{W_{j-1}}{W_j}}\}$.

\bex
Let $(R,\cm)$ be a local ring as in example \ref{Ex.Rings.with.Taylor.approximation}.
\\1. Let $G=\mathcal{R}$ be the group of local coordinate changes, $\bx\to \phi(\bx)$. They act on the elements of the ring by $f(\bx)\to\phi^*(f(\bx))=f(\phi(\bx))$.
For the filtration $\{\cm^j\}$ the group  $G^0$ consists of the changes of the form $\bx\to\bx+h(\bx)$, where $h(\bx)\in\cm^2$.
\\2. More generally, consider the group of automorphisms of a module, $G=GL(p,R)\rtimes \mathcal{R}\circlearrowright R^{\oplus p}$, acting by $w(\bx)\to U(\bx)w(\phi(x))$, where $\phi\in \mathcal{R}$, while $U(\bx)$ is an invertible matrix over $R$. Then
\begin{equation*}
G^0=\{(U,\phi),\ \phi(\bx)=\bx+h(\bx),\ h\in\cm^2,\ U(\bx)=\one+u(\bx),\ u(\bx)\in Mat(p,p;\cm)\}.
\end{equation*}
3. Note that $G^0$ depends essentially on the filtration. In the previous examples we could take the filtration by the powers of some other ideal, $\{J^i\}$, or just by a decreasing sequence of ideals.
\eex

\subsubsection{Logarithm, exponent and the tangent space}\label{Sec.Group.Action.Logarithm.Exponent.Tangent.Space}
As is mentioned after theorem \ref{Thm.Determinacy} we can pass to the completion of the module, $\hW$, \wrt the $\{W_j\}$ filtration. Accordingly we have $GL(\hW_\bullet)\supset\hG\supseteq\hG^0$, the completions of
$GL(W)$, $G$, $G^0$.

Among all the $\k$-linear maps (not necessarily invertible), $End_\k(\hW)$, consider the nilpotent ones,  $End^{nilp}_\k(\hW):=\{\xi\big| \ \xi \hW_j\sseteq \hW_{j+1}\}$. Consider the logarithmic map (recall that $\Q\sset \k$):
\begin{equation*}
\hG^0\stackrel{ln}{\to}End^{nilp}_\k(\hW),\quad
g\to ln(g):=\suml^\infty_{k=1}\frac{(1-g)^k}{k}.
\end{equation*}
As $g\in\hG^0$, $(1-g)\hW_j\sseteq \hW_{j+1}$, thus the sum (though infinite) is a well defined linear operator on $\hW$.
As this logarithm is defined by the standard formula, we have $ln(a^ia^j)=ln(a^i)+ln(a^j)$. But in general $ln(ab)\neq ln(a)+ln(b)$, as $a,b$ do not commute. Nevertheless we assume:
\beq\label{Eq.Assumption.Group.Has.Logarithm}
\text{the image $ln(\hG^0)$ is a $\k$-linear subspace of $End^{nilp}_\k(\hW)$.}
\eeq
This is satisfied in many cases, e.g. in all our examples.
\bed
The tangent space to $\hG^0$ at the unit element is the $\k$-module $T_{\hG^0}:=ln(\hG^0)\sseteq  End_\k(\hW)$.
\eed
Consider the exponential map, $T_{\hG^0}\stackrel{exp}{\to}GL^0(W_\bullet)$, defined by
 $exp(\xi):=\one+\suml^\infty_{k=1}\frac{\xi^k}{k!}$. As $\xi$ is a nilpotent endomorphism, the sum
  (though infinite) is a well defined linear operator and is invertible.
\bel\label{Thm.Group.Action.Exponent.Log.Inverse}
$exp\Big(T_{\hG^0}\Big)=\hG^0$ and the maps $T_{\hG^0}\underset{ln}{\stackrel{exp}{\rightleftarrows}}\hG^0$ are the mutual inverses.
\eel
\bpr
Let $\xi\in T_{\hG^0}$, then $\xi=ln(g)$, for some $g\in\hG^0$. Thus $\exp(\xi)=exp(ln(g))=g\in \hG^0$.

The maps $ln$ and $exp$ are mutual inverses as they are defined by the same Taylor series as the classical functions.
\epr

\subsubsection{The relevant properties of the exponent and the variation operator}
The $j$'th stabilizer of $w\in W$ is the subgroup $St_j(w)=\{g\in G: \pi_j(gw)=\pi_j(w)\}$.
For any $g\in G^0$ and $w\in W$ define the variation operator, $\De_w(g):=gw-w$.
\bel\label{Thm.Linearization.Lemma.1}
The restriction $\pi_{j+1}\De_w|_{St_j(w)\cap G^0}\to\pi_{j+1}W_j$, $j\ge1$ is a homomorphism of groups.
\eel
\bpr
First note that the image of $\De_w|_{St_j(w)}$ is indeed in $W_j$, as $\pi_j(\De_w(g))=0$ for any $g\in St_j(w)$.

Let $g\in G^0$ and $h\in St_j(w)$. Then
\begin{equation*}
\pi_{j+1}\De_w(gh)=\pi_{j+1}\big(\De_{hw}(g)+\De_w(h)\big)=\pi_{j+1}\big(\De_{w}(g)\big)+\pi_{j+1}\big(\De_w(h)\big),
\end{equation*}
 as $\pi_j(hw)=\pi_j(w)$.
\epr

\bel\label{Thm.Linearization.Lemma.2} Let $\xi\in T_{\hG^0}$ and $w\in W$.
\\1.  $\pi_jexp(\xi)\in \pi_j(St_j(w))$ iff $\pi_j(\xi w)=0\in\pi_j(W)$.
\\2. If $\pi_j exp(\xi)\in \pi_{j} St_{j}(w)$ then $\pi_{j+1}\De_w(exp(\xi))=\pi_{j+1}(\xi w)$.
\eel
\bpr
1. $\Rrightarrow$
As the stabilizer is a group, $\pi_jexp(t\xi)w=\pi_j w$ for all $t\in\Z$. The left hand side of this equation is a polynomial in $t$ because $\xi$ is nilpotent. As $char(\k)=0$ and $\k\supseteq\Q$, the equality holds for all $t\in\k$. But this implies $\pi_j \xi w=0$.

$\Lleftarrow$ If $\pi_j\xi w=0$ then $\pi_j\xi^k w=0$, thus $\pi_j exp(\xi)\in\pi_j St_j(w)$.

2.The function $h(t)=\pi_{j}\De_w(exp(t\xi))$ is polynomial in $t$. By lemma \ref{Thm.Linearization.Lemma.1} it is additive. Thus $h(t)=tc$ where $c=h(1)=\pi_j (exp(\xi w)-w)=\pi_{j}(\xi w)$.
\epr

\subsection{Proof of theorem \ref{Thm.Determinacy}}\label{Sec.Group.Action.Proof.of.Theorem}
As is explained after the statement of the theorem, it is enough to consider the completion, $\hG^0\circlearrowright\hW$. Thus we can use the exponent map $T_{\hG^0}\stackrel{exp}{\to}\hG^0$, cf. \S\ref{Sec.Group.Action.Logarithm.Exponent.Tangent.Space}.
  Fix some $w\in \hW$ and consider the corresponding maps $T_{\hG^0}\stackrel{L(\la)=\la(w)}{\underset{F(\la)=exp(\la)(w)-w}{\rightrightarrows}}\hW$.
  Note that $L(T_{\hG^0})=T_{(\hG^0w,w)}$ and $F(T_{\hG^0})=\hG^0w-w$. Then the theorem can be formulated in the form:
\begin{equation*}
\text{If $\hW_k\sseteq L(T_{\hG^0})$ then $\hW_k\sseteq F(T_{\hG^0})$.}
\end{equation*}
Note that $F(0)=0\in \hW$. The statement will follow from lemma \ref{Thm.Group.Action.IFT} if we show:
\begin{equation*}
\text{$F(T_{\hG^0})+\Big(L(T_{\hG^0})\cap \hW_j\Big)\sseteq F(T_{\hG^0})+\hW_{j+1}$ for any $j\ge k$.}
\end{equation*}
 Thus we should check that for any $\mu,\la_j\in T_{\hG^0}$ such that 
  $L(\la_j)=\la(w_j)\in \hW_j$ there exists $\mu'\in T_{\hG^0}$ such that
\begin{equation*}
 \Big(exp(\mu)(w)-w\Big)+\la_j(w)\in \Big(exp(\mu')(w)-w\Big)+\hW_{j+1}.
\end{equation*}
 Define $\mu'$ by $exp(\mu')=exp(\la_j)exp(\mu)$, by lemma \ref{Thm.Group.Action.Exponent.Log.Inverse} such $\mu'$ exists and is unique. Note that $\pi_j(\la_j(w))=0$ for the chosen $w\in\hW$.
 Then, by lemma \ref{Thm.Linearization.Lemma.1}: \begin{equation*}
 \pi_{j+1}(exp(\la_j)exp(\mu)(w)-w)=\pi_{j+1}(exp(\la_j)(w)-w)+\pi_{j+1}(exp(\mu)(w)-w).
 \end{equation*}
   Further, as $\la_j(w)\in \widehat{W_j}$ we get by lemma \ref{Thm.Linearization.Lemma.2}: $\pi_{j+1}(exp(\la_j)(w)-w)=\pi_{j+1}(\la_j(w))$.
 Altogether
\begin{equation*}
 exp(\la_j)exp(\mu)(w)-w+\widehat{W_{j+1}}=exp(\mu)(w)-w+\la_j(w)+\widehat{W_{j+1}},
 \end{equation*}
  as needed.
\epr

\subsection{An application to finite determinacy}\label{Sec.Group.Action.Determinacy.Appplications}
Let $W$ be a module over a local ring $(R,\cm)$, with the filtration $\{W_j=\cm^j W\}$. Suppose $G$ preserves the filtration, then we get:
\begin{equation*}
\text{ If $\cm^k W\sseteq\overline{T_{(G^0w,w)}}$ then $w+\cm^kW\sseteq\overline{G^0w}$.}
\end{equation*}
When $R,W$ are complete  this is a ready criterion, otherwise one uses the Artin approximation theorem (or theorem \ref{Thm.Approximation.Theorem.for.smooth.funcs} in the $C^\infty$-case). This recovers the classical criterion of \cite{Mather.I}-\cite{Mather.III}, revised and generalized many times (\cite{Gaffney}, \cite{Wall}, \cite{Damon.1984}): {\em the determinacy is fixed on the tangent level}.

\

Let $(R,\cm)$ be as in example \ref{Ex.Rings.with.Taylor.approximation}. Below we describe several scenarios (the module and the group action), in each case it is enough to write down the corresponding tangent space(s).
\bex
  Let $W=R$, so one studies the determinacy of function germs. The group of local coordinate changes, $\bx\to\bx+\phi(\bx)$, acts by $f(\bx)\to f(\bx+\phi(\bx))$ and induces the right equivalence, $\cR$. The contact equivalence, $\cK$, is induced by $f(\bx)\to (1+u(\bx))f(\bx+\phi(\bx))$. The unipotent parts, $\cR^0$, $\cK^0$, are  realized for $u(\bx)\in\cm$, $\phi(\bx)\in\cm^2$. Denote by $Der(R,\cm^2)$ the $R$-module of all derivations from $R$ to $\cm^2$. The tangent space to the orbit is then $T_{(\cK^0f,f)}=Der(R,\cm^2)(f)+\cm(f)$. Thus we get, compare e.g. to \cite[Theorem I.2.23]{GLS}:
\begin{equation*}\ber
\text{If $\cm^k\sseteq Der(R,\cm^2)(f)$ then $f$ is $k$-$\cR^0$-determined.}\\
\text{If $\cm^k\sseteq Der(R,\cm^2)(f)+\cm(f)$ then $f$ is $k$-$\cK^0$-determined.}
\eer\end{equation*}
\eex
\bex
More generally, let $W=R^{\oplus p}$ with the filtration $\cm^j W$. The contact group action can be written as $f\to (\one+U)f(\bx+\phi(\bx))$, where $U\in GL(p,R)$, $\phi\in Aut(R)$. For the unipotent part, $\cK^0$, one has: $U\in Mat(p,p;\cm)$, $\phi\in\cm^2$. Then $T_{(\cK^0f,f)}=Der(R,\cm^2)f+Mat(p,p;\cm)f$.
\eex
\bex
Consider the $R$-module of matrices, $Mat(m,n;R)$. The groups $G_{lr}=GL(m,R)\times GL(n,R)$, $\cR$, $\cG_{lr}:=G_{lr}\rtimes \mathcal{R}$ act on $Mat(m,n;R)$ by $A(\bx)\to UA(\phi(\bx))V^{-1}$. Then
\begin{equation*}
T_{(\cG^0_{lr}A,A)}=Span_R\Big(UA-AV\Big)_{(U,V)\in Mat(m,m;\cm)\times Mat(n,n;\cm)}+Span_R\Big(\cD(A)\Big)_{\cD\in Der(R,\cm^2)}.
\end{equation*}
This group and various its subgroups are important in many areas. The determinacy questions are studied in \cite{BK.I}, \cite{BK.II}.
\eex
\bex
When the hypersurface singularity $\{f=0\}\sset Spec(R)$ is non-isolated, the tangent space $T_{(\cK^0f,f)}$ does not contain $\cm^k$ for any $k$. Thus the filtration $\{\cm^j\}$ is irrelevant.  It is natural to consider only the deformations preserving the singular locus. More precisely, for the ideal $Jac_f+(f)$ consider the following saturation. Take the primary decomposition $\capl_i I_i$ and apply the procedure: if $\sqrt{I_i}\supsetneq\sqrt{I_j}$ then erase $I_i$ in this decomposition. Eventually one gets a saturated version, $(Jac_f+(f))^{sat}$, geometrically this corresponds to removing the embedded components of lower dimension. Then one can consider either of the filtrations:
 $W_j=\cm^{j-1}(Jac_f+(f))^{sat}$ or $W_j=\cm^{j-1}\int (Jac_f+(f))^{sat}$,
 here $\int I=\{g\in R\big| \ Der(R)(g)+(g)\sseteq I\}$. The later filtration has been studied in \cite{Siersma1983},\cite{Siersma1987}, \cite{Pelikaan}.

In both cases one defines $Der_{W_1}(R):=\{\cD\in Der(R,\cm^2)\big| \ \cD W_1\sseteq\cm W_1\}$ and considers the corresponding subgroup $exp(Der_{W_1}(R))$ of $\cR^0$.  In both cases one has:
\begin{equation*}
\text{If $Der_{W_1}(R)(f)\supseteq\cm^k W_1$ then $f$ is $k$-determined for deformations inside $W_1$.}
\end{equation*}
\eex

\end{document}